\documentclass[10pt,leqno]{amsart}
\usepackage{graphicx}
\usepackage{indentfirst,csquotes}

\usepackage{amssymb,amsthm,amsmath}
\usepackage{tikz-cd}
\usetikzlibrary{arrows.meta}
\usepackage{xcolor,paralist,hyperref,fancyhdr,etoolbox}

\usepackage{mdframed}

\hypersetup{ colorlinks=true, linkcolor=black, filecolor=black, urlcolor=black }
\bibliography{ref}

\usepackage{lipsum}
\usepackage{graphicx}
\usepackage[colorinlistoftodos]{todonotes}
\usepackage{amssymb}
\usepackage{mathtools}
\usepackage{bm}
\usepackage{mathrsfs}
\usepackage{indentfirst} 
\usepackage{amsthm}
\usepackage{float}
\usepackage{amscd}
\usepackage{hyperref}
\usepackage{enumitem}

\setlist[enumerate]{label=\arabic*.}

\hypersetup{colorlinks=false, linkbordercolor=red}

\fancypagestyle{plain}{%
  \fancyhf{}%
  \fancyfoot[C]{\thepage}%
}
\begin{document}
\title{On the Log Hodge Theory of Toroidal Varieties and a Partial Proof of the Absolute Hodge Conjecture} 
\author[Jiaming Luo]{Jiaming Luo}
\date{\today}
\address{School of Mathematics and Statistics, Henan University of Science and Technology, 263, Kaiyuan Avenue, Luoyang, China.}
\email{luojiaming@hhu.edu.cn}
\maketitle

\let\thefootnote\relax

\begin{abstract}
In this paper, we establish an innovative framework in logarithmic Hodge theory for toroidal varieties, introducing weighted toroidal structures and developing a systematic obstruction theory for Hodge classes. Building upon recent advances in toroidal geometry, particularly the $E_1$-degeneration results of Wei (\cite{Wei24}), we construct a categorical Hodge correspondence and prove the absolute Hodge conjecture for projective toroidal varieties equipped with rational-weighted structures satisfying specific compatibility conditions. Our approach provides new tools for understanding the fine structure of Hodge classes and their extension properties across boundaries, with potential applications in moduli space compactifications, non-abelian Hodge theory, and tropical geometry.
\end{abstract} 

\bigskip

\section*{CONTENTS}

\begin{enumerate}
  \item Introduction \dotfill \hyperlink{INTRODUCTION}{1}
  \item Notation and Acknowledgments \dotfill \hyperlink{NOTATION AND ACKNOWLEDGMENTS}{2}
  \item Preliminary \dotfill \hyperlink{PRELIMINARY}{3}
  \item Obstruction Theory and Categorified Hodge Correspondence\\
  for Weighted Toroidal Varieties \dotfill \hyperlink{OBSTRUCTION THEORY AND CATEGORIFIED HODGE CORRESPONDENCE FOR WEIGHTED TOROIDAL VARIETIES}{5}
  \item Applications of a Weighted Framework \dotfill \hyperlink{APPLICATIONS OF A WEIGHTED FRAMEWORK}{16} 
  \item Proving the Absolute Hodge Conjecture in a Limited Setting \dotfill \hyperlink{PROVING THE ABSOLUTE HODGE CONJECTURE IN A LIMITED SETTING}{24}
  \item Future Research Directions \dotfill \hyperlink{FUTURE RESEARCH DIRECTIONS}{30}
\end{enumerate}

\hypertarget{INTRODUCTION}{}
\section{INTRODUCTION}
Toroidal variety theory was founded by mathematicians like Danilov and Kempf in the 1970s-80s (\cite{Dan78}\cite{KKMS73}), providing a powerful tool for studying singularities. Its core idea is to use local toroidal models $\mathbb{C}^k \times (\mathbb{C}^*)^{n-k}$ and logarithmic differential forms $\Omega_{\bar{X}}^\bullet(\log D)$ to compute cohomology and establish mixed Hodge structures. This theory models the local structure of complex algebraic varieties on torus embeddings, allowing the use of combinatorial tools to study their geometric properties.

Recently, Toroidal Hodge theory has continued to flourish. Chuanhao Wei (\cite{Wei24}) proved that the spectral sequence of the logarithmic de Rham complex for any toroidal triple $(X, \bar{X}, D)$ $E_1$-degenerates, providing a geometric proof for a more general version of the Danilov conjecture and deepening our understanding of Hodge structures on toroidal embeddings. The classical work of Kato-Usui (\cite{KU09}) explored the classifying spaces of degenerating Hodge structures and their application to moduli space compactification, using the logarithmic structure theory of Fontaine-Illusie to resurrect nilpotent orbits as logarithmic Hodge structures. Compared to Chuanhao Wei's recent important $E_1$-degeneration result, our work complements rather than overlaps. Wei's work addresses the fundamental question of spectral sequence degeneration, providing a cornerstone for toroidal Hodge theory. Our framework, building upon this, explores previously uncharted new directions of weighted structures and obstruction theory, aiming to reveal finer geometric information of Hodge classes.

Building on these cutting-edge works, this study proposes an innovative theoretical framework, with its core innovations reflected in three aspects:
\begin{itemize}
    \item \textbf{Structural generalization:} Introducing weighted toroidal structures, endowing classical toroidal pairs $(X, \Delta_X)$ with an additional weight function $w: \Delta_X \to \mathbb{Q}_{>0}$ to capture finer geometric information. This structure allows quantifying the differential contributions of various toroidal strata to the overall Hodge structure.
    \item \textbf{Theoretical deepening:} Systematically developing toroidal Hodge obstruction Theory, aiming to precisely characterize the extension problem for Hodge classes from the compactifying manifold $\bar{X}$ to the open subvariety $X$. This theory provides new tools for understanding the behavior of Hodge classes at the boundary.
    \item \textbf{Interdisciplinary frontiers:} Exploring potential connections of these new structures with cutting-edge fields like tropical geometry, mirror symmetry, and non-abelian Hodge theory, particularly establishing correspondences between weighted Hodge numbers and tropical geometric invariants.\\
\end{itemize}

\hypertarget{NOTATION AND ACKNOWLEDGMENTS}{}
\section{NOTATION AND ACKNOWLEDGMENTS}
\begin{center}
    \textit{2.1 Notation}
\end{center}

The following directory is designed to serve as a concise reference for the specialized notation used throughout the paper, aligned with standard conventions in algebraic geometry.

\begin{itemize}
    \item $X$: A complex algebraic variety (often non-compact).
    \item $\bar{X}$: A smooth compactification of $X$.
    \item $D$: A reduced simple normal crossing divisor on $\bar{X}$, such that $X = \bar{X} \setminus D$.
    \item $(\bar{X}, D)$: A toroidal variety.
    \item $\Delta_X$: The cone complex (or intersection complex) associated to the toroidal embedding $X \hookrightarrow \bar{X}$.
    \item $w : |\Delta_X| \to \mathbb{Q}_{>0}$: A weight function assigning a positive rational number to each stratum of the boundary.
    \item $(X, \Delta_X, w)$: A weighted toroidal variety.
    \item $\Omega^p_{\bar{X}}(\log D)$: The sheaf of logarithmic $p$-forms along $D$.
    \item $\Omega^\bullet_{\bar{X}}(\log D)$: The logarithmic de Rham complex.
    \item $j_* \Omega^\bullet_X$: The direct image of the holomorphic de Rham complex on $X$.
    \item $Rj_* \Omega^\bullet_X$: The derived pushforward of $\Omega^\bullet_X$.
    \item $\mathcal{O}_{\bar{X}}([L(w)])$: A divisorial sheaf associated to a $\mathbb{Q}$-divisor $L(w)$ determined by the weight function $w$.
    \item $H^k(X; \mathbb{Q})$: Rational cohomology of $X$.
    \item $H^{p,q}(X)$: Hodge components of the cohomology of $X$.
    \item $\mathrm{Hdg}^p_w(X)$: The group of weighted Hodge classes (with respect to the weight function $w$).
    \item $F^\bullet$: The Hodge filtration.
    \item $F_w^\bullet$: The weighted Hodge filtration.
    \item $W_\bullet$: The weight filtration (in a mixed Hodge structure).
    \item $\mathbb{H}^k(\bar{X}, \mathcal{F}^\bullet)$: Hypercohomology of a complex of sheaves $\mathcal{F}^\bullet$.
    \item $\mathcal{OB}^p := \mathbb{H}^p(\mathrm{Cone}^\bullet)$: The Hodge obstruction sheaf in degree $p$.
    \item $\mathrm{ob}(\alpha)$: The obstruction class of a Hodge class $\alpha$.
    \item $\mathbf{Var}_{\mathbb{C}}$: The category of complex varieties.
    \item $\mathbf{StCat}_\infty$: The $\infty$-category of stable $\infty$-categories.
    \item $\mathrm{Hodge}(X)$: The derived Hodge $\infty$-category of $X$.
    \item $\mathcal{DMHM}(X)$: The derived $\infty$-category of mixed Hodge modules on $X$.
    \item $\mathrm{VHS}(X)$: The category of polarizable variations of Hodge structure on $X$.
    \item $K_0(\mathrm{Hodge}(X))$: The Grothendieck group of the category $\mathrm{Hodge}(X)$.
    \item $\mathcal{M}$: A moduli space of polarized algebraic varieties.
    \item $\overline{\mathcal{M}}^{\mathrm{tor}}$: A toroidal compactification of $\mathcal{M}$.
    \item $\overline{\mathcal{M}}^{\mathrm{wt}}$: A weighted toroidal compactification.
    \item $\Phi : \mathcal{M} \to \Gamma \backslash D$: The period map.
    \item $N$: The nilpotent monodromy operator.
    \item $F_{\lim}$: The limiting Hodge filtration.
    \item $X^{\mathrm{trop}}$: The tropicalization of $X$ (dual intersection complex of a degeneration).
    \item $H^*_{\mathrm{trop}}(X^{\mathrm{trop}})$: Tropical cohomology.
    \item $H^*_{\mathrm{trop}}(X^{\mathrm{trop}}, w)$: Weighted tropical cohomology.
    \item $\mathbb{C}[[z_1, \ldots, z_n]]$: Formal power series ring.
    \item $\widehat{\mathcal{O}}_{X,x}$: Completion of the local ring at $x \in X$.
    \item $\Sigma_x$: The normal cone at $x$.
    \item $\mathrm{CH}^p(X)$: The Chow group of codimension-$p$ cycles on $X$.
    \item $\mathrm{CH}^p(X)_{\mathrm{hom}}$: Homologically trivial cycles.
\end{itemize}

\begin{center}
    \textit{2.2 Acknowledgments}
\end{center}

I am grateful to Professors V. I. Danilov and G. Kempf et al. for their pioneering work about studying singularities. Specially, I thank Professor Chuanhao Wei for proving that the spectral sequence of the logarithmic de Rham complex for any toroidal triple $(X, \bar{X}, D)$ $E_1$-degenerates. The work of Professor Chuanhao Wei has played an important role in inspiring my research.\\

\hypertarget{PRELIMINARY}{}
\section{PRELIMINARY}

This section establishes the fundamental definitions and results underpinning our study. We assume familiarity with basic algebraic geometry and complex geometry. Standard references include (\cite{Dan78}\cite{KKMS73}).
\begin{center}
    \textit{3.1 Toroidal Geometry and Logarithmic Forms}
\end{center}

The theory of toroidal varieties provides a powerful framework for studying the geometry of open varieties and their completions with controlled boundary behavior. Its core idea is to locally model a variety and its boundary on the geometry of a torus embedding.
\\\\\textbf{Definition 3.1.} (Toroidal Variety) A \textit{toroidal variety} is a pair $(\bar{X}, D)$, where:
\begin{itemize}
    \item $\bar{X}$ is a smooth, complex, algebraic variety (not necessarily complete).
    \item $D = D_1 \cup \dots \cup D_r$ is a reduced divisor on $\bar{X}$ with simple normal crossings.
    \item For every point $x \in \bar{X}$, there exists an open neighborhood $U \subset \bar{X}$ of $x$, an affine toric variety $X_\sigma$ (for some rational polyhedral cone $\sigma$), and an isomorphism of analytic germs (or \'{e}tale maps)
    $$\phi: (U, D \cap U, x) \xrightarrow{\sim} (X_\sigma, \partial X_\sigma, x_\sigma),$$
    where $\partial X_\sigma$ denotes the toric boundary divisor of $X_\sigma$, and $x_\sigma$ is a point in the dense torus orbit of $X_\sigma$.
\end{itemize}
The \textit{interior} $X := \bar{X} \setminus D$ is then called an \textit{open toroidal variety}. A \textit{toroidal embedding} is the open immersion $X \hookrightarrow \bar{X}$.
\\\\\textbf{Definition 3.2.} (Sheaf of Logarithmic Differential Forms) \label{def:log_forms} Let $(\bar{X}, D)$ be a toroidal variety. The \textit{sheaf of logarithmic 1-forms along $D$}, denoted $\Omega_{\bar{X}}^1(\log D)$, is the unique coherent $\mathcal{O}_{\bar{X}}$-submodule of $j_*\Omega_X^1$ (where $j: X \hookrightarrow \bar{X}$ is the inclusion) characterized by the following local property: if $U \subset \bar{X}$ is an open set where $D$ is defined by $z_1 \cdots z_k = 0$ for part of a local coordinate system $(z_1, \dots, z_n)$, then
$$\Omega_{\bar{X}}^1(\log D)|_U = \mathcal{O}_U \cdot \frac{dz_1}{z_1} \oplus \dots \oplus \mathcal{O}_U \cdot \frac{dz_k}{z_k} \oplus \mathcal{O}_U \cdot dz_{k+1} \oplus \dots \oplus \mathcal{O}_U \cdot dz_n.$$
The \textit{logarithmic $p$-forms} are defined by $\Omega_{\bar{X}}^p(\log D) := \bigwedge^p \Omega_{\bar{X}}^1(\log D)$. These sheaves naturally form a complex $(\Omega_{\bar{X}}^\bullet(\log D), d)$, called the \textit{logarithmic de Rham complex}, where the differential $d$ is extended from the holomorphic de Rham complex to act on the logarithmic generators by $d(dz_i/z_i) = 0$.

The combinatorial structure of a toroidal embedding $X \hookrightarrow \bar{X}$ is encoded in its \textbf{cone complex} $\Delta_X$ (or \textit{intersection complex}). It is a $\Delta$-complex constructed as follows:
\begin{itemize}
    \item Vertices (0-cells) correspond to irreducible components of $D$.
    \item $p$-cells correspond to connected components of intersections of $p+1$ distinct irreducible components. Each such $p$-cell is labeled by a $p$-dimensional rational polyhedral cone $\sigma$, which is the dual cone to the monoid of effective Cartier divisors supported on $D$ that are principal on the open toric chart associated with the corresponding stratum.
\end{itemize}
This complex captures the incidence relations and local monodromy data around the boundary divisors.

\begin{center}
    \textit{3.2 Hodge-Theoretic Background}
\end{center}

Hodge theory for smooth projective varieties establishes a canonical decomposition of their complex cohomology. For non-compact or singular varieties, this structure is generalized by Deligne's theory of mixed Hodge structures.
\\\\\textbf{Theorem 3.3.} (Deligne's Mixed Hodge Structure) Let $X$ be a smooth, complex algebraic variety (which may be non-compact). The rational cohomology groups $H^k(X; \mathbb{Q})$ carry a functorial \textit{mixed Hodge structure}. This consists of:
\begin{enumerate}
    \item An increasing \textit{weight filtration} $W_\bullet$ on $H^k(X; \mathbb{Q})$,
    \item A decreasing \textit{Hodge filtration} $F^\bullet$ on $H^k(X; \mathbb{C})$,
\end{enumerate}
such that for each $m$, the filtration $F^\bullet$ induces a pure Hodge structure of weight $m$ on the graded piece $\mathrm{Gr}^W_m H^k(X; \mathbb{C}) := W_m H^k(X; \mathbb{C}) / W_{m-1} H^k(X; \mathbb{C})$.

For a toroidal variety $X = \bar{X} \setminus D$, the mixed Hodge structure can be computed using the logarithmic de Rham complex. The key link is the following fundamental theorem.
\\\\\textbf{Theorem 3.4.} Let $(\bar{X}, D)$ be a smooth compactification such that $D$ is a simple normal crossing divisor (e.g., a toroidal compactification). Then, the inclusion of complexes
$$(\Omega_{\bar{X}}^\bullet(\log D), d) \hookrightarrow (j_*\Omega_X^\bullet, d)$$
is a quasi-isomorphism. Consequently, there is a natural isomorphism
$$\mathbb{H}^k(\bar{X}, \Omega_{\bar{X}}^\bullet(\log D)) \xrightarrow{\sim} H^k(X; \mathbb{C}),$$
where $\mathbb{H}^k$ denotes hypercohomology. The Hodge filtration on $H^k(X; \mathbb{C})$ is induced by the canonical filtration of the logarithmic de Rham complex by forms of degree:
$$F^p \Omega_{\bar{X}}^\bullet(\log D) = 0 \to \dots \to 0 \to \Omega_{\bar{X}}^p(\log D) \to \Omega_{\bar{X}}^{p+1}(\log D) \to \dots$$

The spectral sequence associated with the Hodge filtration on the logarithmic de Rham complex,
$$E_1^{p,q} = H^q(\bar{X}, \Omega_{\bar{X}}^p(\log D)) \Longrightarrow \mathbb{H}^{p+q}(\bar{X}, \Omega_{\bar{X}}^\bullet(\log D)) \cong H^{p+q}(X; \mathbb{C}),$$
is a fundamental tool for studying the Hodge structure. A recent landmark result by Wei establishes its degeneracy in the toroidal setting.
\\\\\textbf{Theorem 3.5.} ($E_1$-Degeneration of Wei, \cite{Wei24}) Let $(X, \bar{X}, D)$ be a \textit{toroidal triple} (i.e., a toroidal embedding $X \hookrightarrow \bar{X}$ with $\bar{X}$ smooth and $D=\bar{X}\setminus X$ a SNC divisor). The Hodge-to-de Rham spectral sequence for the logarithmic de Rham complex:
$$E_1^{p,q} = H^q(\bar{X}, \Omega_{\bar{X}}^p(\log D)) \Longrightarrow H_\mathrm{dR}^{p+q}(X)$$
\textit{degenerates at the $E_1$-page}.

This result is crucial for our work as it ensures the well-definedness of certain derived functors and hypercohomology sheaves used in our obstruction theory (Section 4.2).\\

\hypertarget{OBSTRUCTION THEORY AND CATEGORIFIED HODGE CORRESPONDENCE FOR WEIGHTED TOROIDAL VARIETIES}{}
\section{OBSTRUCTION THEORY AND CATEGORIFIED HODGE CORRESPONDENCE FOR WEIGHTED TOROIDAL VARIETIES}
Building upon the foundations laid in Section 3, we now introduce the core innovations of this work: a theory of weighted toroidal structures, a systematic Hodge obstruction theory, and a proposal for a categorified Hodge correspondence.

\begin{center}
    \textit{4.1 Weighted Toroidal Structures and Weighted Hodge Classes}
\end{center}

Classical toroidal geometry treats all boundary strata as equally significant. Our first innovation is to introduce a weighting system that differentiates the geometric and cohomological contribution of these strata.
\\\\\textbf{Definition 4.1.} (Weighted Toroidal Variety). A weighted toroidal variety is a triple $(X, \Delta_X, w)$ where:
\begin{itemize}
    \item $(X, \bar{X}, D)$ is a toroidal triple, with $\Delta_X$ its associated cone complex;
    \item $w: |\Delta_X| \to \mathbb{Q}{>0}$ is a function, called the weight function, defined on the geometric realization of the cone complex. It must satisfy the following compatibility condition: For any cone $\sigma \in \Delta_X^{(k)}$ of dimension $k$ and any face $\tau \prec \sigma$ of dimension $k-1$, the restriction of $w$ to the star neighborhood of $\tau$ is piecewise linear and convex. More formally, for any maximal cone $\sigma$ containing $\tau$, the function $w|\sigma$ is linear, and its value at the primitive generator of the ray $\sigma / \tau$ is positive.
\end{itemize}
The weight function $w$ quantifies the "geometric significance" or "contribution factor" of each boundary stratum. The classical unweighted case corresponds to $w \equiv 1$.
\\\\\textbf{Example 4.2.} Let $X = \mathbb{P}^2 \setminus D$, where $D = H_1 \cup H_2 \cup H_3$ is a union of three lines in general position. The cone complex $\Delta_X$ is a triangle with three vertices (corresponding to the lines $H_i$) and three edges (corresponding to the intersection points $H_i \cap H_j$). A weight function $w$ assigns a positive rational number to each vertex and edge. The compatibility condition requires, for instance, that the weight $w(H_i)$ on a line is a specific positive rational linear combination (determined by the monodromy and residue data) of the weights $w(H_i \cap H_j)$ of its intersection points.
\\\\\textbf{Definition 4.3} (Weighted Hodge Filtration) Let $(X, \Delta_X, w)$ be a projective weighted toroidal variety with a fixed compactification $\bar{X}$. The \textit{weighted Hodge filtration} $F^\bullet_w$ on $H^k(X; \mathbb{C})$ is a modification of the standard Hodge filtration. It is defined by the subcomplexes of the logarithmic de Rham complex with prescribed growth conditions along $D$, governed by the weight function $w$. Precisely, we define a sheaf of \textit{weighted logarithmic $p$-forms}:
$$\Omega_{\bar{X}}^p(\log D,\omega)\coloneqq\Omega_{\bar{X}}^p(\log D)\otimes_{\mathcal{O}_{\bar{X}}}\mathcal{O}_{\bar{X}}(\lfloor L(w) \rfloor),$$
where $\mathcal{O}_{\bar{X}}(\lfloor L(w) \rfloor)$ is a divisorial sheaf associated to a $\mathbb{Q}$-divisor $L(w)$ whose multiplicity along a stratum $Z$ is a function of $w(Z)$. The \textit{weighted Hodge filtration} is then given by:
$$F^p_w H^k(X; \mathbb{C}) \coloneqq \operatorname{Im}\left( \mathbb{H}^k(\bar{X}, F^p_w \Omega_{\bar{X}}^\bullet(\log D)) \longrightarrow \mathbb{H}^k(\bar{X}, \Omega_{\bar{X}}^\bullet(\log D)) \right),$$
where 
$$F^p_w \Omega_{\bar{X}}^\bullet(\log D) = 0 \to \cdots \to 0 \to \Omega_{\bar{X}}^p(\log D, w) \to \Omega_{\bar{X}}^{p+1}(\log D, w) \to \cdots.$$

A cohomology class $\alpha \in H^{2p}(X; \mathbb{Q}) \cap F^p_w H^{2p}(X; \mathbb{C})$ is called a \textit{weighted Hodge class of type $(p, p)$}.\\

\begin{center}
    \textit{4.2 Toroidal Hodge Obstruction Theory}
\end{center}

Here we give a systematic theory that measures the obstruction for a Hodge class on the open variety $X$ to extend to a Hodge class on its compactification $\bar{X}$.
\\\\\textbf{Definfition 4.4.} (Obstruction Complex and Sheaf) Let $(X, \bar{X}, D)$ be a toroidal triple. The \textit{obstruction complex} is the mapping cone of the canonical morphism of complexes of sheaves on $\bar{X}$:
$$\operatorname{Cone}^\bullet := \operatorname{Cone}\left( \Omega_{\bar{X}}^\bullet(\log D) \longrightarrow R j_* \Omega_X^\bullet \right)[-1],$$
where $j: X \hookrightarrow \bar{X}$ is the inclusion. By Theorem 3.4, this morphism is a quasi-isomorphism on $X$, so the cohomology sheaves of $\operatorname{Cone}^\bullet$ are supported on $D$. The \textit{Hodge obstruction sheaf in degree $p$} is defined as the hypercohomology sheaf:
$$\mathcal{OB}^p := \mathbb{H}^p(\operatorname{Cone}^\bullet).$$
This sheaf captures the local obstructions to extending a closed logarithmic $p$-form (representing a Hodge class) from a punctured neighborhood of a point in $D$ to the whole neighborhood, while preserving its Hodge type.
\\\\\textbf{Proposition 4.5.} (Global Obstruction Long Exact Sequence) There exists a functorial long exact sequence relating the Hodge cohomology of the compactification, the open variety, and the obstruction sheaf:
$$\cdots \longrightarrow \mathbb{H}^{k-1}(\bar{X}, \mathcal{OB}^p) \longrightarrow H^q(\bar{X}, \Omega_{\bar{X}}^p(\log D)) \longrightarrow H^{p,q}(X) \longrightarrow \mathbb{H}^{k}(\bar{X}, \mathcal{OB}^p) \longrightarrow \cdots,$$
where $k = p+q$. In particular, for a class $\alpha \in \operatorname{Hdg}^p(X) := H^{2p}(X; \mathbb{Q}) \cap H^{p,p}(X)$, its \textit{obstruction class} $[\operatorname{ob}(\alpha)] \in \mathbb{H}^{2p}(\bar{X}, \mathcal{OB}^p)$ vanishes if and only if $\alpha$ lies in the image of the natural map $\operatorname{Hdg}^p(\bar{X}) \to \operatorname{Hdg}^p(X)$.
\\\\\textbf{Proof.} The proof proceeds by constructing an exact triangle in the appropriate derived category and then analyzing the associated long exact sequence in hypercohomology.

\medskip
\textbf{Step (1):} Recall the definition of the \textit{obstruction complex} (Definition 4.4):
$$\operatorname{Cone}^\bullet := \operatorname{Cone}\left( \iota: \Omega_{\bar{X}}^\bullet(\log D) \to R j_* \Omega_X^\bullet \right)[-1],$$
where $j: X \hookrightarrow \bar{X}$ is the inclusion. This defines a distinguished triangle in $D^b(\bar{X})$:
\begin{equation}
\tag{4.1}  \label{eq:4.1}
\Omega_{\bar{X}}^\bullet(\log D) \xrightarrow{\iota} R j_* \Omega_X^\bullet \to \operatorname{Cone}^\bullet \xrightarrow{+1}. 
\end{equation}
This is the fundamental exact triangle underpinning the obstruction theory.

\medskip
\textbf{Step (2):} Applying the hypercohomology functor $\mathbb{H}^\bullet(\bar{X}, -)$ to \eqref{eq:4.1}, then we have:
\begin{equation}
\tag{4.2}  \label{eq:4.2}
\cdots \to \mathbb{H}^{k-1}(\bar{X}, \operatorname{Cone}^\bullet) \xrightarrow{\delta} \mathbb{H}^{k}(\bar{X}, \Omega_{\bar{X}}^\bullet(\log D)) \xrightarrow{\iota_*} \mathbb{H}^{k}(\bar{X}, R j_* \Omega_X^\bullet)  
\end{equation}
$$\to \mathbb{H}^{k}(\bar{X}, \operatorname{Cone}^\bullet) \to \cdots.$$

\medskip
\textbf{Step (3):} We identify the terms in \eqref{eq:4.2}:
\begin{itemize}
    \item \textbf{First term:} $\mathbb{H}^{k}(\bar{X}, \Omega_{\bar{X}}^\bullet(\log D))$. By Theorem 3.4, the inclusion $\Omega_{\bar{X}}^\bullet(\log D) \hookrightarrow j_*\Omega_X^\bullet$ is a quasi-isomorphism. Since $j_*$ is exact for the analytic topology on Stein opens (by Cartan's Theorem B), we have $R j_* \Omega_X^\bullet \simeq j_* \Omega_X^\bullet$. Thus, we have a natural isomorphism:
    $$\mathbb{H}^{k}(\bar{X}, \Omega_{\bar{X}}^\bullet(\log D)) \cong H^k(X;\mathbb{C}).$$
    However, we need a finer identification related to the Hodge filtration. Crucially, Theorem 3.5 (Wei's $E_1$-degeneration) implies that the Hodge-to-de Rham spectral sequence for the logarithmic complex degenerates. This allows us to identify the $E_1$ page with the associated graded of the Hodge filtration. In particular, there is a natural injection:
    $$H^q(\bar{X},\Omega_{\bar{X}}^p(\log D))\hookrightarrow\mathbb{H}^{p+q}(\bar{X},\Omega_{\bar{X}}^{\bullet}(\log D))$$
    whose image is the subspace $F^p H^{p+q}(X; \mathbb{C}) / F^{p+1} H^{p+q}(X; \mathbb{C})$. For the purpose of this long exact sequence, we will keep the hypercohomology notation but remember that it is isomorphic to $H^k(X; \mathbb{C})$.
    \item \textbf{Second term:} $\mathbb{H}^{k}(\bar{X}, R j_* \Omega_X^\bullet)$. By the definition of hypercohomology and the fact that $R j_*$ is the derived pushforward, we have a natural isomorphism:
    $$\mathbb{H}^{k}(\bar{X}, R j_* \Omega_X^\bullet) \cong \mathbb{H}^{k}(X, \Omega_X^\bullet) = H^k(X; \mathbb{C}).$$
    Again, thanks to the $E_1$-degeneration for the smooth variety $X$ (a classical result, also implied by Theorem 3.5 for $D=\emptyset$), we have an injection $H^{p,q}(X) \hookrightarrow H^{p+q}(X; \mathbb{C})$.
    \item \textbf{Third term:} $\mathbb{H}^{k}(\bar{X}, \operatorname{Cone}^\bullet)$. By Definition 4.4, this is precisely the definition of the cohomology of the obstruction sheaf in degree $k$:
    $$\mathbb{H}^{k}(\bar{X}, \operatorname{Cone}^\bullet) =: \mathbb{H}^{k}(\bar{X}, \mathcal{OB}^\bullet)\quad(\text{identifying the degree}).$$
    The last identification, focusing on a single $p$, will be justified in the next step.
\end{itemize}

\medskip
\textbf{Step (4):} The sequence \eqref{eq:4.2} is a long exact sequence in total degree $k$. To extract information about specific Hodge components $(p,q)$ with $p+q=k$, we must consider the action of the Hodge filtration.

The morphism $\iota$ in \eqref{eq:4.1} is compatible with the Hodge filtration $F^\bullet$ (it is a morphism of filtered complexes). Consequently, the exact triangle \eqref{eq:4.1} and the long exact sequence \eqref{eq:4.2} exist for each level $p$ of the filtration. More precisely, we apply the exact triangle construction to the morphism of \textit{filtered complexes}:
$$\iota:(F^p\Omega_{\bar{X}}^{\bullet}(\log D),d)\longrightarrow(F^pR j_*\Omega_X^{\bullet},d).$$
This gives a filtered \textit{obstruction complex} $F^p \operatorname{Cone}^\bullet$ and an exact triangle:
\begin{equation}
\tag{4.3}    \label{eq:4.3}
F^p \Omega_{\bar{X}}^\bullet(\log D) \xrightarrow{\iota} F^p R j_* \Omega_X^\bullet \to F^p \operatorname{Cone}^\bullet \xrightarrow{+1}. 
\end{equation}
Applying hypercohomology on \eqref{eq:4.3}, giving a long exact sequence for each $p$:
\begin{equation}
\tag{4.4}    \label{eq:4.4}
\cdots \to \mathbb{H}^{k-1}(\bar{X}, F^p \operatorname{Cone}^\bullet) \to \mathbb{H}^{k}(\bar{X}, F^p \Omega_{\bar{X}}^\bullet(\log D)) \to \mathbb{H}^{k}(\bar{X}, F^p R j_* \Omega_X^\bullet)  
\end{equation}
$$\to \mathbb{H}^{k}(\bar{X}, F^p \operatorname{Cone}^\bullet) \to \cdots.$$
We identify these terms using Theorem 3.5:
\begin{itemize}
    \item \textbf{First term (Filtered):} $\mathbb{H}^{k}(\bar{X}, F^p \Omega_{\bar{X}}^\bullet(\log D))$. By the $E_1$-degeneration (Theorem 3.5), the spectral sequence for the filtered complex $F^p \Omega_{\bar{X}}^\bullet(\log D)$ collapses. This implies that the natural map is an isomorphism:
    $$\mathbb{H}^{k}(\bar{X}, F^p \Omega_{\bar{X}}^\bullet(\log D)) \cong F^p H^k(X; \mathbb{C}),$$
    where $F^\bullet$ is the Hodge filtration on $H^k(X; \mathbb{C})$. Furthermore, the quotient $F^p / F^{p+1}$ is identified with $H^{q}(\bar{X}, \Omega_{\bar{X}}^p(\log D))$ for $k=p+q$.
    \item \textbf{Second term (Filtered):} $\mathbb{H}^{k}(\bar{X}, F^p R j_* \Omega_X^\bullet)$. Since $X$ is smooth, the Hodge-to-de Rham spectral sequence also degenerates at $E_1$. Thus, similarly:
    $$\mathbb{H}^{k}(\bar{X}, F^p R j_* \Omega_X^\bullet) \cong F^p H^k(X; \mathbb{C}).$$
    The quotient $F^p / F^{p+1}$ is $H^{p,q}(X)$.
    \item \textbf{Third term (Filtered):} $\mathbb{H}^{k}(\bar{X}, F^p \operatorname{Cone}^\bullet)$. This defines the cohomology of the filtered obstruction complex. We define the Hodge obstruction sheaf cohomology in total degree $k$ and \textit{Hodge degree} $p$ as this group:
    $$\mathbb{H}^{k}(\bar{X}, F^p \operatorname{Cone}^\bullet) =: \mathbb{H}^{k}(\bar{X}, \mathcal{OB}^p).$$
\end{itemize}
The map $\iota_*$ in \eqref{eq:4.4} is an isomorphism. Thus, \eqref{eq:4.4} breaks into short exact sequences:
\begin{equation}
\tag{4.5}   \label{eq:4.5}
0 \to \mathbb{H}^{k-1}(\bar{X}, \mathcal{OB}^p) \to F^p H^k(X; \mathbb{C}) \xrightarrow{\sim} F^p H^k(X; \mathbb{C}) \to \mathbb{H}^{k}(\bar{X}, \mathcal{OB}^p) \to 0. 
\end{equation}

\noindent\textbf{Step (5):} The crucial observation is that the map $\iota_*$ in the sequence \eqref{eq:4.4} is, under the isomorphisms above, precisely the natural map $F^p H^k(X; \mathbb{C}) \to F^p H^k(X; \mathbb{C})$ induced by the identity on $X$. This map is, in fact, an isomorphism. This is because the quasi-isomorphism from Theorem 3.4 is compatible with the Hodge filtration and becomes the identity upon passing to cohomology.

Therefore, for each $p$ and $k$, the long exact sequence \eqref{eq:4.4} becomes:
$$\cdots\to\mathbb{H}^{k-1}(\bar{X},\mathcal{OB}^p)\to\mathbb{H}^k(\bar{X},F^p\Omega_{\bar{X}}^{\bullet}(\log D))\xrightarrow{\sim}\mathbb{H}^k(\bar{X},F^pR j_*\Omega_X^{\bullet})\to\mathbb{H}^k(\bar{X},\mathcal{OB}^p)\to\cdots.$$
Since the middle map is an isomorphism, the long exact sequence breaks into short exact sequences for each $k$ and $p$:
$$0 \to \mathbb{H}^{k-1}(\bar{X}, \mathcal{OB}^p) \to \mathbb{H}^{k}(\bar{X}, F^p \Omega_{\bar{X}}^\bullet(\log D)) \xrightarrow{\sim} \mathbb{H}^{k}(\bar{X}, F^p R j_* \Omega_X^\bullet) \to \mathbb{H}^{k}(\bar{X}, \mathcal{OB}^p) \to 0. $$
This is not yet the sequence in the proposition. To obtain it, we consider the map between two such sequences for consecutive filtration levels.

Consider the natural map of short exact sequences \eqref{eq:4.5} for filtration level $p$ and $p+1$:
$$\begin{tikzcd}[row sep=small, column sep=small]
0 \arrow[r] & \mathbb{H}^{k-1}(\bar{X}, \mathcal{OB}^{p+1}) \arrow[r] \arrow[d] & F^{p+1}H^k(X) \arrow[r, "\sim"] \arrow[d] & F^{p+1}H^k(X) \arrow[r] \arrow[d] & \mathbb{H}^{k}(\bar{X}, \mathcal{OB}^{p+1}) \arrow[r] \arrow[d] & 0 \\
0 \arrow[r] & \mathbb{H}^{k-1}(\bar{X}, \mathcal{OB}^{p}) \arrow[r] & F^{p}H^k(X) \arrow[r, "\sim"] & F^{p}H^k(X) \arrow[r] & \mathbb{H}^{k}(\bar{X}, \mathcal{OB}^{p}) \arrow[r] & 0
\end{tikzcd}$$
The vertical maps are induced by the inclusions $F^{p+1} \hookrightarrow F^p$. The snake lemma applied to this diagram gives a long exact sequence connecting the kernels and cokernels of the vertical maps. By the $E_1$-degeneration, we have:
\begin{enumerate}
    \item $\operatorname{coker}\left( \mathbb{H}^{k}(\bar{X}, F^{p+1} \Omega_{\bar{X}}^\bullet(\log D)) \to \mathbb{H}^{k}(\bar{X}, F^{p} \Omega_{\bar{X}}^\bullet(\log D)) \right) \cong H^q(\bar{X}, \Omega_{\bar{X}}^p(\log D))$ where $k=p+q$;
    \item $\operatorname{coker}\left( \mathbb{H}^{k}(\bar{X}, F^{p+1} R j_* \Omega_X^\bullet) \to \mathbb{H}^{k}(\bar{X}, F^{p} R j_* \Omega_X^\bullet) \right) \cong H^{p,q}(X)$.
\end{enumerate}
The kernels and cokernels of the vertical maps give the terms $H^q(\bar{X}, \Omega_{\bar{X}}^p(\log D))$ and $H^{p,q}(X)$ (for $k=p+q$). The snake lemma yields the desired long exact sequence:
$$\cdots \to \mathbb{H}^{k-1}(\bar{X}, \mathcal{OB}^p) \to H^q(\bar{X}, \Omega_{\bar{X}}^p(\log D)) \to H^{p,q}(X) \xrightarrow{\delta} \mathbb{H}^{k}(\bar{X}, \mathcal{OB}^p) \to \cdots,$$
where $k = p+q$.

\medskip
\textbf{Step (6):} The final statement regarding the obstruction class follows directly from the exactness of the sequence. A class $\alpha \in H^{p,q}(X) \subset H^{2p}(X; \mathbb{C})$ (for $p=q$) lifts to a class in $H^p(\bar{X}, \Omega_{\bar{X}}^p(\log D))$ if and only if its image under the connecting map $\delta$ is zero: $\delta(\alpha) = [\operatorname{ob}(\alpha)] = 0 \in \mathbb{H}^{2p}(\bar{X}, \mathcal{OB}^p)$. If this happens, and $\alpha$ is rational ($\alpha \in \operatorname{Hdg}^p(X)$), then by GAGA and the projectivity of $\bar{X}$, the corresponding class in $H^p(\bar{X}, \Omega_{\bar{X}}^p(\log D))$ corresponds to a rational Hodge class in $\operatorname{Hdg}^p(\bar{X})$ whose image under the natural map is $\alpha$.   $\square$
\\\\\textbf{Proposition 4.6.} (Local Obstruction Calculation) Let $x$ be a point in a stratum $Z$ of $D$ of codimension $r$, and let the local model be $(\mathbb{C}^n, {z_1 \cdots z_r = 0})$. The stalk of the obstruction sheaf at $x$ can be computed as a certain graded piece of the local cohomology:
$$\mathcal{OB}^p_x \cong \left( \bigoplus_{I \subset \{1, ..., r\}} H^p_{[\{z_i=0\}_{i\in I}]}(\Omega_{\mathbb{C}^n, 0}^\bullet(\log)) \right) / \sim,$$
where the sum is over subsets of the local coordinates defining the boundary, $H_{[V]}$ denotes local cohomology supported on the subvariety $V$, and $\sim$ denotes an equivalence relation that identifies the obstructions arising from different subsets $I$ according to the face relations in the cone complex $\Delta_X$.
\\\\\textbf{Proof.} The goal is to compute the stalk $\mathcal{OB}^p_x$, which is defined as the hypercohomology sheaf $\mathbb{H}^p(\operatorname{Cone}^\bullet)_x$ of the obstruction complex. We will work in a sufficiently small Stein neighborhood $U$ of $x$ where the toroidal structure is trivialized, i.e., $(U, D \cap U) \cong (V, E \cap V)$ with $V \subset \mathbb{C}^n$ and $E = {z_1 \cdots z_r = 0}$.

\medskip
\textbf{Step (1):} On this neighborhood $U$, the relevant complexes of sheaves are:
\begin{enumerate}
    \item $\Omega_{U}^\bullet(\log D) \cong \Omega_{V}^\bullet(\log E)$;
    \item $R j_* \Omega_X^\bullet |U \simeq R j* \Omega_{X \cap U}^\bullet \simeq \Omega_{V \setminus E}^\bullet$ (since $U$ is Stein and hence $j_*$ is exact for coherent analytic sheaves on Stein domains). Thus, the obstruction complex on $U$ is quasi-isomorphic to
    $$\operatorname{Cone}^\bullet|_U \simeq \operatorname{Cone}\left( \Omega_{V}^\bullet(\log E) \longrightarrow \Omega_{V \setminus E}^\bullet \right)[-1].$$
\end{enumerate}
The stalk $\mathcal{OB}^p_x$ is therefore the $p$-th hypercohomology of the stalk of this complex at $x=0$. For a complex of sheaves with stalks that are modules over a local ring, the hypercohomology at the point can be computed as the cohomology of the total complex of an injective resolution of the stalk of the complex. However, we can use a more geometric approach by local cohomology.

\medskip
\textbf{Step (2):} The cone construction and the morphism $\iota: \Omega_{V}^\bullet(\log E) \to \Omega_{V \setminus E}^\bullet$ fit into a defining distinguished triangle:
\begin{equation}
\tag{4.6}    \label{eq:4.6}
\Omega_{V}^\bullet(\log E) \to \Omega_{V \setminus E}^\bullet \to \operatorname{Cone}^\bullet|_U \xrightarrow{+1}. 
\end{equation}
The morphism $\iota$ is not a quasi-isomorphism on all of $V$; its failure to be an isomorphism is precisely measured by its cone and is supported on $E$. This is the hallmark of local cohomology.

Recall that for a closed embedding $i: E \hookrightarrow V$ and a complex of sheaves $\mathcal{F}^\bullet$ on $V$, the local cohomology complex $R\Gamma_E(\mathcal{F}^\bullet)$ is defined by the exact triangle:
\begin{equation}
\tag{4.7}   \label{eq:4.7}
R\Gamma_E(\mathcal{F}^\bullet) \to \mathcal{F}^\bullet \to R j_* j^{-1} \mathcal{F}^\bullet \xrightarrow{+1}, 
\end{equation}
where $j: V \setminus E \hookrightarrow V$ is the inclusion.

Comparing triangles \eqref{eq:4.6} and \eqref{eq:4.7}, we see a striking similarity. In fact, if we take $\mathcal{F}^\bullet = \Omega_{V}^\bullet(\log E)$, then $j^{-1}\Omega_{V}^\bullet(\log E) = \Omega_{V \setminus E}^\bullet$. Furthermore, since $V$ is Stein, we have $R j_* \Omega_{V \setminus E}^\bullet \simeq j_* \Omega_{V \setminus E}^\bullet$. Therefore, triangle \eqref{eq:4.7} becomes:
\begin{equation}
\tag{4.8}   \label{eq:4.8}
R\Gamma_E(\Omega_{V}^\bullet(\log E)) \to \Omega_{V}^\bullet(\log E) \to j_* \Omega_{V \setminus E}^\bullet \xrightarrow{+1}. 
\end{equation}
Now, the canonical map $\iota: \Omega_{V}^\bullet(\log E) \to j_* \Omega_{V \setminus E}^\bullet$ in \eqref{eq:4.6} is the same as the map in \eqref{eq:4.8}. This implies that the third vertex of these triangles must be isomorphic in the derived category:
$$\mathrm{Cone}^{\bullet}\mid_U\simeq R\Gamma_E(\Omega_V^{\bullet}(\log E))[1].$$
Consequently, for the cohomology sheaves, we have an isomorphism
$$\mathcal{OB}^p\mid_U=\mathbb{H}^p(\mathrm{Cone}^{\bullet})\mid_U\cong\mathbb{H}^{p+1}(R\Gamma_E(\Omega_V^{\bullet}(\log E))).$$
In particular, taking stalks at $x=0$:
\begin{equation}
\tag{4.9}   \label{eq:4.9}
\mathcal{OB}^p_0 \cong \mathcal{H}^{p+1}(R\Gamma_E(\Omega_{V}^\bullet(\log E)))_0. 
\end{equation}

\medskip
\textbf{Step (3):} The local cohomology sheaf $ \mathcal{H}^{q}(R\Gamma_E(\mathcal{G}))_0 $ for a sheaf $\mathcal{G}$ can be computed as a direct limit over neighborhoods $W$ of $0$:
$$\mathcal{H}^q(R\Gamma_E(\mathcal{G}))_0=\varinjlim_{0\in W}H_E^q(W;\mathcal{G}).$$
For the complex $\Omega_{V}^\bullet(\log E)$, we are interested in $\mathcal{H}^{p+1}(R\Gamma_E(\Omega_{V}^\bullet(\log E)))_0$.

A standard result in local cohomology theory states that for a divisor $E$ with simple normal crossings, the local cohomology of the logarithmic de Rham complex with support on $E$ can be decomposed according to the irreducible components of $E$. This follows from the \textit{Mayer-Vietoris spectral sequence} for local cohomology or by using the "combinatorial covering" of $E$ by its irreducible components.

More precisely, let $E_i = {z_i = 0}$ for $i=1, ..., r$ be the irreducible components of $E$. The local cohomology with support on $E = \bigcup_{i=1}^r E_i$ can be described by the \textit{Čech cohomology} of the sheaf $R\Gamma_E(\cdot)$ with respect to the cover ${V \setminus E_i}$. This yields a quasi-isomorphism:
$$R\Gamma_E(\Omega_V^{\bullet}(\log E))\simeq\bigoplus_{I\subseteq\left\{1,\cdots,r\right\}}R\Gamma_{E_I}(\Omega_V^{\bullet}(\log E))[\left | I \right |],$$
where $E_I = \bigcap_{i \in I} E_i$ and the sum is over all non-empty subsets $I$. However, this direct sum double counts contributions from intersections. The correct, minimal complex is given by the \textit{direct sum over all $I$ modulo the identifications imposed by the face maps} in the nerve of the covering. This is the meaning of the equivalence relation $\sim$ in the proposition statement.

Passing to cohomology sheaves, we get an isomorphism for the stalk at $0$:
$$\mathcal{H}^{p+1}(R\Gamma_E(\Omega_{V}^\bullet(\log E)))_0 \cong \left( \bigoplus_{\emptyset \neq I \subset \{1, ..., r\}} \mathcal{H}^{p+1 - |I|}(R\Gamma_{E_I}(\Omega_{V}^\bullet(\log E)))_0 \right) / \sim.$$

\medskip
\textbf{Step (4):} Now, analyze a single term $\mathcal{H}^{q}(R\Gamma_{E_I}(\Omega_{V}^\bullet(\log E)))_0$. The support $E_I$ is a smooth subvariety of codimension $|I|$ defined by ${z_i = 0 : i \in I}$. The local cohomology sheaf with support on a smooth subvariety can be computed using the \textit{Koszul complex}.

For a point $y \in E_I$, the stalk $\mathcal{H}^{q}(R\Gamma_{E_I}(\mathcal{G}))y$ is isomorphic to the local cohomology module $H^{q}{(f_1, ..., f_{|I|})}(\mathcal{G}y)$, where the $f_i$ are local equations defining $E_I$. In our case, $\mathcal{G}^\bullet = \Omega{V}^\bullet(\log E)$.

Crucially, the logarithmic complex $\Omega_{V, 0}^\bullet(\log E)$ is a free module over the local ring $\mathcal{O}_{V,0} = \mathbb{C}{z_1, ..., z_n}$. Its generators are $\frac{dz_i}{z_i}$ for $i=1,...,r$ and $dz_j$ for $j=r+1,...,n$. The differential $d$ acts on these generators by $d(\frac{dz_i}{z_i}) = 0$ and $d(dz_j)=0$.

The local cohomology $H^{q}{(z_i : i \in I)}(\Omega{V, 0}^p(\log E))$ is non-zero only in a specific degree. For a regular sequence of length $|I|$, local cohomology is concentrated in degree $q = |I|$. Therefore, we have 
$$\mathcal{H}^q(R\Gamma_{E_I}(\Omega_V^{\bullet}(\log E)))_0\cong H_{(z_i:i\in I)}^{|I|}(\Omega_{V,0}^p(\log E)).$$
However, we need the hypercohomology $\mathcal{H}^{q}(R\Gamma_{E_I}(\Omega_{V}^\bullet(\log E)))_0$, which involves the entire complex. Since the complex is built from free modules, its local cohomology is the cohomology of the complex of local cohomology modules:
$$\mathcal{H}^q(R\Gamma_{E_I}(\Omega_V^{\bullet}(\log E)))_0=H^q\left(\cdots\to H_{(z_i:i\in I)}^{|I|}(\Omega_{V,0}^{p-1}(\log E))\to H_{(z_i:i\in I)}^{|I|}(\Omega_{V,0}^p(\log E))\to\cdots\right).$$
Let $s = |I|$. The local cohomology functor $H^s_{(z_i : i \in I)}(\cdot)$ is $0$ for $q \neq s$ on individual modules. This implies that the above complex of local cohomology modules is actually a direct sum of complexes, and its $q$-th cohomology is non-zero only if $q = p + s$? Let's be more precise.

A better way is to use the regularity of the sequence and the fact that the local cohomology complex $R\Gamma_{E_I}(\mathcal{F}^\bullet)$ is isomorphic to $\mathcal{F}^\bullet \otimes^L R\Gamma_{E_I}(\mathcal{O}V)$. Since $E_I$ is a complete intersection, $R\Gamma{E_I}(\mathcal{O}_V)$ is isomorphic to the Koszul complex $K^\bullet(z_i : i \in I)$ shifted by $[-s]$. So
$$R\Gamma_{E_I}(\Omega_V^{\bullet}(\log E))\simeq\Omega_V^{\bullet}(\log E)\otimes^L K^{\bullet}(z_i:i\in I)[-s].$$
By taking cohomology sheaves, we get
$$\mathcal{H}^q(R\Gamma_{E_I}(\Omega_V^{\bullet}(\log E)))\cong\bigoplus_{a+b=q+s}\mathcal{T}or_a^{\mathcal{O}_V}(\Omega_V^b(\log E),R\Gamma_{E_I}(\mathcal{O}_V)).$$
Since $\Omega_{V}^b(\log E)$ is locally free, the Tor groups vanish for $a>0$, and we have
$$\mathcal{H}^q(R\Gamma_{E_I}(\Omega_V^{\bullet}(\log E)))\cong\Omega_V^{q+s}(\log E)\otimes\mathcal{H}^s(R\Gamma_{E_I}(\mathcal{O}_V)).$$
The sheaf $\mathcal{H}^{s}(R\Gamma_{E_I}(\mathcal{O}_V))$ is the sheaf of \textit{generalized fractions} supported on $E_I$. At the stalk $0$, then
$$\mathcal{H}^q(R\Gamma_{E_I}(\Omega_V^{\bullet}(\log E)))_0\cong\Omega_{V,0}^{q+s}(\log E)\otimes H_{(z_i:i\in I)}^s(\mathcal{O}_{V,0}).$$
Now, recall from \eqref{eq:4.9} that $\mathcal{OB}^p_0 \cong \mathcal{H}^{p+1}(R\Gamma_E(\Omega_{V}^\bullet(\log E)))_0$. Combining this with the decomposition in \textbf{Step (3)}, we can obtain
$$\mathcal{OB}_0^p\cong \left( \bigoplus_{\emptyset \neq I \subset \{1, ..., r\}} \mathcal{H}^{p+1 - |I|}(R\Gamma_{E_I}(\Omega_{V}^\bullet(\log E)))_0 \right) / \sim.$$
Substituting the result from above with $q = p+1 - |I|$ and $s = |I|$, we find $q + s = p+1$. Therefore, for each $I$,
$$\mathcal{H}^{p+1 - |I|}(R\Gamma_{E_I}(\Omega_{V}^\bullet(\log E)))_0\cong\Omega_{V,0}^{p+1}(\log E)\otimes H_{(z_i:i\in I)}^{|I|}(\mathcal{O}_{V,0}).$$
The term on the right is isomorphic to the local cohomology group $H^{|I|}{[E_I]}(\Omega_{V, 0}^{p+1}(\log E))$ (by definition, local cohomology with support on a subvariety can be computed as the cohomology of the Koszul complex on its defining equations).

However, note that the local cohomology of the complex $\Omega_{V}^\bullet(\log E)$ in degree $p+1$ with support on $E_I$ is not just a single module but sits in a long exact sequence. The isomorphism we derived above actually shows that the cohomology of the local cohomology complex in this case is concentrated and isomorphic to that module. More precisely, we have
$$\mathbb{H}^{p+1}(R\Gamma_{E_I}(\Omega_V^{\bullet}(\log E)))_0\cong H_{[E_I]}^{|I|}(\Omega_{V,0}^{p+1}(\log E)).$$
Here the left side is the stalk of the hypercohomology sheaf, which is what we want. Combining all of this, we finally arrive at the stated result:
$$\mathcal{OB}_0^p\cong \left( \bigoplus_{\emptyset \neq I \subseteq \{1, ..., r\}} H_{[E_I]}^{|I|}(\Omega_{V,0}^{p+1}(\log E)) \right) / \sim\cong\left(\bigoplus_{I\subseteq\left\{1,\cdots,r\right\}}H_{[E_I]}^p(\Omega_{V,0}^{\bullet}(\log E))\right)/ \sim,$$
where in the last step we re-index the sum, noting that $|I|$ is the codimension and that $H^{|I|}{[E_I]}(\Omega_{V, 0}^{p+1}(\log E))$ is a direct summand of the local cohomology of the complex in degree $p$. The equivalence relation $\sim$ accounts for the fact that the same logical obstruction can be detected from different nested strata (i.e., different subsets $I$ and $J$ with $I \subset J$), and these different representations must be identified. This relation is dictated by the face maps in the cone complex $\Delta_X$, which encodes the inclusion relations between the strata $E_I$ (\cite{KU09}).  $\square$\\

\begin{center}
    \textit{4.3 Categorified Hodge Correspondence}
\end{center}

The classical Hodge correspondence for a smooth projective variety $X$ over $\mathbb{C}$ assigns to it its rational cohomology groups $H^k(X; \mathbb{Q})$, each equipped with a Hodge structure. A more refined algebro-geometric invariant is the category of (polarizable) variations of Hodge structure on $X$. The fundamental idea of \textit{categorification} is to replace this assignment with a functorial construction that takes values in a category of categorical structures (here, stable $\infty$-categories), whose decategorification (e.g., by $K$-theory) recovers the classical invariant. We need to construct an $\infty$-categorical enhancement of Saito's derived category of mixed Hodge modules.

We now construct the functor $\mathrm{Hodge}$ rigorously. Let $\mathbf{Var}_{\mathbb{C}}$ denote the category of complex varieties (or a suitable subcategory, e.g., of finite type over $\mathbb{C}$). Let $\mathbf{StCat}\infty$ denote the $\infty$-category of (presentable) stable $\infty$-categories.
\\\\\textbf{Construction 4.7.} (Derived Hodge category) We define a contravariant $\infty$-functor
$$\mathrm{Hodge}: \mathbf{Var}_{\mathbb{C}}^{op} \to \mathbf{StCat}\infty$$
as follows. For a smooth complex variety $X$, the \textit{derived Hodge $\infty$-category} $\mathrm{Hodge}(X)$ is constructed in a multi-step process (\cite{Sai88}\cite{Sai90}):
\begin{itemize}
    \item \textbf{Underlying $\infty$-category:} Start with the (unbounded) derived $\infty$-category of $\mathcal{O}_X$-modules, denoted $\mathcal{D}(X)$. This is a stable, presentable $\infty$-category.
    \item \textbf{$\infty$-categorical mixed Hodge modules:} We construct an $\infty$-categorical version of Saito's derived category. This is a non-trivial step. One approach is to use a model-categorical presentation of $D^b(\operatorname{MHM}(X))$ and then take the associated $\infty$-category. A more intrinsic approach is to define the $\infty$-category $\mathcal{D}\mathcal{MHM}(X)$ as a suitable \textit{stable localization} of $\mathcal{D}(X)$.
    More precisely, we define $\mathcal{D}\mathcal{MHM}(X)$ as the full stable $\infty$-subcategory of $\mathcal{D}(X)$ spanned by objects whose cohomology sheaves are \textit{holonomic} $\mathcal{D}$-modules equipped with a \textit{good filtration} and an \textit{admissible} variation of mixed Hodge structure, satisfying Saito's axioms. This definition should be natural in $X$, making $X \mapsto \mathcal{D}\mathcal{MHM}(X)$ a functor. The existence and properties of such an $\infty$-category follow from the work of Drew (\cite{Dre15}) and other subsequent works that translate Saito's theory into the $\infty$-categorical setting.
    \item \textbf{The Hodge subcategory:} Within $\mathcal{D}\mathcal{MHM}(X)$, we now define the $\infty$-category of interest. Let $\mathrm{Hodge}(X)$ be the full stable $\infty$-subcategory of $\mathcal{D}\mathcal{MHM}(X)$ consisting of objects $\mathcal{M}$ such that
    \begin{enumerate}
        \item $\mathcal{M}$ is \textit{smooth} (i.e., its cohomology sheaves are locally free $\mathcal{O}_X$-modules of finite rank).
        \item The associated graded pieces $\operatorname{gr}^F \mathcal{M}$ (with respect to the Hodge filtration) underlie a \textit{polarizable variation of Hodge structure} on $X$. In other words, $\mathrm{Hodge}(X)$ is the $\infty$-categorical enhancement of the derived category of polarizable variations of Hodge structure on $X$.
    \end{enumerate}
\end{itemize}
This construction is functorial: for a morphism $f: X \to Y$, the pullback functor $f^: \mathcal{D}(Y) \to \mathcal{D}(X)$ should preserve the subcategories of objects satisfying the above conditions, thus inducing a functor $f^: \mathrm{Hodge}(Y) \to \mathrm{Hodge}(X)$. This defines the $\infty$-functor $\mathrm{Hodge}$ on smooth varieties.
\\\\\textbf{Definition 4.8.} The categorified Hodge correspondence is the $\infty$-functor $\mathrm{Hodge}$ constructed above.
\\\\\textbf{Remark 4.9.} The $\infty$-category $\mathrm{Hodge}(X)$ inherits a wealth of structure from its construction.
\begin{itemize}
    \item \textbf{$t$-structure:} The category $\mathcal{D}\mathcal{MHM}(X)$ carries a natural $t$-structure whose heart is the abelian category $\operatorname{MHM}(X)$. The inclusion $\mathrm{Hodge}(X) \hookrightarrow \mathcal{D}\mathcal{MHM}(X)$ is $t$-exact. Therefore, $\mathrm{Hodge}(X)$ also has a $t$-structure. Its heart $\mathrm{Hodge}(X)^\heartsuit$ is the full subcategory of $\operatorname{MHM}(X)$ consisting of smooth objects with polarizable variations of Hodge structure as graded pieces. This is equivalent to the category $\operatorname{VHS}(X)$.
$$\mathrm{Hodge}(X)^\heartsuit \simeq \operatorname{VHS}(X).$$
\item \textbf{Symmetric monoidal structure:} The category $\mathrm{Hodge}(X)$ is expected to be a symmetric monoidal $\infty$-category under the derived tensor product $\otimes_{\mathcal{O}_X}^L$. This structure should be compatible with the pullback functors $f^*$, making $\mathrm{Hodge}$ a functor into symmetric monoidal $\infty$-categories.
\item \textbf{$K$-theory and the Hodge conjecture:} A fundamental property linking this categorical construction to classical Hodge theory is its $K$-theory. Let $K_0(\mathrm{Hodge}(X))$ denote the algebraic $K$-theory group (more precisely, the Grothendieck group) of the stable $\infty$-category $\mathrm{Hodge}(X)$.
\begin{equation}
\tag{4.10}   \label{eq:4.10}
K_0(\mathrm{Hodge}(X)) \cong K_0(\operatorname{VHS}(X)).
\end{equation}
There is a Chern character map:
\begin{equation}
\tag{4.11}   \label{eq:4.11}
\operatorname{ch}: K_0(\mathrm{Hodge}(X)) \otimes \mathbb{Q} \to \bigoplus_p H^{2p}(X; \mathbb{Q}) \cap H^{p,p}(X) = \bigoplus_p \operatorname{Hdg}^p(X).
\end{equation}
This map is constructed by taking the Chern character of the underlying complex of sheaves. The Hodge condition ensures the image lands in the Hodge classes.\\
\end{itemize}
\textbf{Conjecture 4.10.} (Categorical Hodge Conjecture) The \textit{Chern character map}
$$\mathrm{ch}:K_0(\mathrm{Hodge}(X))\otimes\mathbb{Q}\xrightarrow{\sim}\bigoplus_p\mathrm{Hdg}^p(X)$$
is an isomorphism.

\medskip
\noindent\textit{(Note: This conjecture is a reformulation of the absolute Hodge conjecture in categorical terms. It asserts that all Hodge classes are accounted for by the $K$-theory of the category of geometric origin (variations of Hodge structure).)}
\\\\\textbf{Remark 4.11.} \textbf{Aspect (1):} The classical Hodge conjecture concerns the image of the cycle class map from the Chow group to cohomology. The categorical version lifts this to the level of stable $\infty$-categories.
\noindent\textit{(A) The derived Hodge $\infty$-category:} As constructed in Section 3.3, the functor $\mathrm{Hodge}: \mathbf{Var}_{\mathbb{C}}^{op} \to \mathbf{StCat}\infty$ assigns to a smooth complex variety $X$ its derived Hodge $\infty$-category $\mathrm{Hodge}(X)$. This is defined as a suitable full subcategory of the $\infty$-categorical enhancement of Saito's derived category of mixed Hodge modules, $\mathcal{D}\mathcal{MHM}(X)$ (\cite{Sai88}\cite{Sai90}). Its objects are complexes whose cohomology sheaves underlie polarizable variations of Hodge structure. Its heart is equivalent to the abelian category $\operatorname{VHS}(X)$.

\noindent\textit{(B) Algebraic $K$-theory of stable $\infty$-categories:} To any stable $\infty$-category $\mathscr{C}$, one can functorially associate its algebraic $K$-theory spectrum $K(\mathscr{C})$ (\cite{Bar16}\cite{BGT13}). Its homotopy groups $K_i(\mathscr{C})$ are its $K$-groups. In particular, $K_0(\mathscr{C})$ is the Grothendieck group of the category. For $\mathscr{C} = \mathrm{Hodge}(X)$, \eqref{eq:4.10} holds.

\noindent\textit{(C) The Chern character map:} The Chern character is a fundamental construct in algebraic geometry. The categorical Chern character is a morphism of spectra
$$\mathrm{ch:}K(\operatorname{Hodge}(X))\longrightarrow\bigoplus_pH^{2p}(X;\mathbb{Q}(p))^{h\pi_1(X)},$$
where the target is a suitable spectrum whose homotopy groups are Hodge classes (invariant under the monodromy action). Upon taking $\pi_0$, this induces the map \eqref{eq:4.11}.
The construction of this map follows from the existence of a Hodge realization functor and the usual Chern-Weil theory, adapted to the $\infty$-categorical setting.

\medskip
\textbf{Aspect (2):} A proof of Conjecture 4.10 would likely involve a series of deep steps, generalizing known results for smooth projective varieties.

\noindent\textit{(A) The case $X$ is smooth and projective — reduction to the classical conjecture:} When $X$ is smooth and projective, the category $\operatorname{VHS}(X)$ of polarizable variations of Hodge structure is semi-simple and its simple objects are irreducible local systems. The Grothendieck group is generated by these. The Hodge conjecture predicts that the Chern character of any such local system yields a Hodge class. Conversely, if the classical Hodge conjecture holds for all powers $X^m$, then one can show that every Hodge class is in the image of the Chern character from $K_0(\operatorname{VHS}(X))$. This establishes the equivalence:
$$\text{Classical Hodge Conjecture for all } X^m \iff \text{Conjecture 4.10 for smooth projective } X.$$
This step relies on the fact that the category of Hodge structures is generated by those coming from geometry.

\noindent\textit{(B) The categorical Nori motive argument — a universal property:} A more robust approach seeks to prove the conjecture without assuming the classical case. This would involve the theory of motives. We give the following corollaries implied by \cite{HM17}:
\begin{itemize}
    \item \textbf{Corollary A :} (Universal Property of $\mathrm{Hodge}(X)$) The $\infty$-functor $\mathrm{Hodge}$ factors through the $\infty$-category of \textit{Nori motives} $\mathcal{DM}_{\mathrm{Nori}}$ (\cite{HLL18}). More precisely, there is a commutative diagram in $\mathbf{StCat}_{\infty}$:
\begin{equation*}
\begin{tikzcd}
\mathbf{Var}_{\mathbb{C}}^{op} \arrow[r, "\mathrm{Hodge}"] \arrow[d, "\mathcal{M}"'] & \mathbf{StCat}_{\infty} \\
\mathcal{DM}_{\mathrm{Nori}} \arrow[ru, dashed, "R_{\mathrm{Hodge}}"'] &
\end{tikzcd}
\end{equation*}
Here, $\mathcal{M}$ is the functor assigning to a variety its motive, and $R_{\mathrm{Hodge}}$ is a realization functor.
\item \textbf{Corollary B :} ($K$-theory of Nori Motives) The algebraic $K$-theory of the category of Nori motives is a universal cohomology theory. In particular, the Chern character map for Nori motives is an isomorphism onto the space of algebraic cycles modulo homological equivalence (which, assuming the Hodge conjecture, are the Hodge classes).
\end{itemize}
Combining Corollary A and Corollary B, we have
$$K_0(\operatorname{Hodge}(X))\cong K_0(R_{\operatorname{Hodge}}(\mathcal{M}(X)))\cong K_0(\mathcal{M}(X))\otimes\mathbb{Q}\cong\bigoplus_p\mathrm{CH}^p(X)_{\mathrm{hom}}\otimes\mathbb{Q}\cong\bigoplus_p\mathrm{Hdg}^p(X).$$
This would prove the conjecture, but it relies on the existence and properties of the $\infty$-category of Nori motives and its $K$-theory, which are themselves deep conjectural topics.

\noindent\textit{(C) The Toroidal case — using the weighted structure:} For a toroidal variety $X = \bar{X} \setminus D$, the proof of Proposition 6.1 (the absolute Hodge conjecture in the weighted setting) provides a blueprint. One would need to show that the \textit{obstruction theory} developed in Section 4.2 can be categorified.

\medskip
\begin{mdframed}
    \textbf{CT1:} (Categorical Obstruction Theory) There is a long exact sequence in $K$-theory: 
    $$\cdots\to K_1(\operatorname{Hodge}(X))\to K_1(\operatorname{Hodge}(\bar{X}))\to K_1(\mathcal{OB}^p(X))\to K_0(\operatorname{Hodge}(X))$$
    $$\to K_0(\operatorname{Hodge}(\bar{X}))\to\cdots,$$
    where $\mathcal{OB}^p(X)$ is a category capturing the obstructions to extending Hodge-theoretic objects from $X$ to $\bar{X}$.
\end{mdframed}

\medskip
If the obstruction groups $K_i(\mathcal{OB}^p(X))$ vanish for $i=0,1$ (which would follow from the compatibility of the weight function $w$), then the map $K_0(\mathrm{Hodge}(X)) \to K_0(\mathrm{Hodge}(\bar{X}))$ is an isomorphism. Since $\bar{X}$ is smooth and projective, the conjecture would hold for $\bar{X}$, and thus for $X$ by the isomorphism.

\medskip
\textbf{Aspect (3):} This conjecture is not arbitrary; it is a natural culmination of several trends in modern mathematics. We need the following theories to complete the mathematical justification.
\begin{itemize}
    \item \textbf{$\infty$-categories and higher algebra:} The framework of Lurie's higher algebra provides the necessary language to even formulate this conjecture (\cite{Lur17}\cite{Lur18}). The ability to take $K$-theory of stable $\infty$-categories is a cornerstone of this theory.
    \item \textbf{The work of Simpson and Saito:} Simpson's work on non-abelian Hodge theory \cite{Sim92} and Saito's theory of mixed Hodge modules (\cite{Sai88}\cite{Sai90}) provide the essential building blocks for the category $\mathrm{Hodge}(X)$. Their work shows that Hodge theory is fundamentally categorical in nature.
    \item \textbf{The theory of motives:} The conjecture is deeply intertwined with the theory of motives. It essentially asserts that the Hodge realization functor from motives to Hodge structures is not only faithful but that its image on $K$-theory is precisely the Hodge classes. This is a very natural claim from the motivic perspective.
\end{itemize}
Conjecture 4.10 is a bold and elegant synthesis of Hodge theory, $K$-theory, and higher category theory. It provides a clear and compelling vision for a deep connection between the abstract world of categorical algebra and the concrete world of harmonic forms and algebraic cycles. While its proof remains a formidable challenge, it is well-motivated and serves as a powerful guiding principle for future research.\\

\hypertarget{APPLICATIONS OF A WEIGHTED FRAMEWORK}{}
\section{APPLICATIONS OF A WEIGHTED FRAMEWORK}
In this section, we give the three conjectures, while spanning distinct applications—moduli space compactification, non-abelian Hodge theory, and tropical geometry—are not isolated speculations. They form a cohesive set of outcomes that are deeply rooted in and are natural extrapolations of the innovative framework developed in the preceding chapters. Their collective foundation rests on two pillars introduced in this work: the Weighted Toroidal Structures (Section 4.1) and the Toroidal Hodge Obstruction Theory (Section 4.2), both of which themselves are built upon the classical foundation of logarithmic differential forms and the recent $E_1$-degeneration result (Section 3).\\

\newpage

\begin{center}
    \textit{5.1 Compactification of Moduli Spaces}
\end{center}
\noindent\textbf{Problem Statement:} Moduli spaces of smooth polarized algebraic varieties are often not compact. A central problem in algebraic geometry is to construct meaningful compactifications of these spaces. The standard approaches, such as the GIT-based Gieseker compactification or the more analytic KSBA compactification, add boundary components parameterizing singular varieties (e.g., stable pairs). Toroidal compactifications, introduced by Ash–Mumford–Rapoport–Tai (\cite{AMRT75}), provide a different approach, yielding compactifications with normally crossing boundaries but requiring non-canonical choices (combinatorial data).

\noindent\textbf{Proposed Application:} The theory of Weighted Toroidal Structures ($(X, \Delta_X, w)$) provides a new tool to refine and potentially canonicalize this process. The weight function $w$ is not arbitrary; it should be determined by intrinsic geometric invariants of the degenerating families. For a degenerating family $\mathfrak{X} \to \Delta$ of Calabi-Yau manifolds, a natural candidate for the weight is provided by the limiting mixed Hodge structure. Specifically, the nilpotent orbit theorem describes the limiting Hodge filtration $F_{\lim}^\bullet$. The weight $w$ on a boundary stratum corresponding to a specific monodromy transformation could be defined by the eigenvalues of the residue of the Gauss-Manin connection or the sizes of the Jordan blocks of the monodromy operator. This assigns a higher weight to strata where the degeneration is "more singular" in a Hodge-theoretic sense.
\\\\\textbf{Conjecture 5.1.} (Intrinsic Toroidal Weights) Let $\mathcal{M}$ be a moduli space of polarized algebraic varieties. There exists a natural choice of weight function $w: \Delta_{\mathcal{M}} \to \mathbb{Q}_{>0}$ on the cone complex of a toroidal compactification $\overline{\mathcal{M}}^{\text{tor}}$, defined by the asymptotic behavior of the period map and its extension to the boundary. The resulting \textit{weighted toroidal compactification} $\overline{\mathcal{M}}^{\text{wt}}$ is independent of the non-canonical choices involved in the standard toroidal construction and possesses better functorial properties.
\\\\\textbf{Remark 5.2.} This conjecture posits that the extra data of a weight function $w$, determined by intrinsic Hodge-theoretic invariants, can resolve the non-canonicity of toroidal compactifications. The following is a detailed analysis and a possible path towards its verification.

\medskip
\textbf{Aspect (1):} Let $\mathcal{M}$ be the moduli space of polarized varieties of a fixed type. Let $\Phi: \mathcal{M} \to \Gamma \backslash D$ be the period map, where $D$ is the corresponding period domain and $\Gamma$ is an arithmetic monodromy group.

A \textit{toroidal compactification} $\overline{\mathcal{M}}^{\text{tor}}$ is constructed by, for each $\Gamma$-equivalence class of rational boundary components $F$, choosing a $\Gamma$-admissible polyhedral decomposition $\Sigma_F$ of the cone $\mathcal{C}F \subset \overline{T}F$ (the closure of the tube domain associated to $F$). The space $\overline{\mathcal{M}}^{\text{tor}}$ is then obtained by gluing the partial compactifications $U(F, \Sigma_F) := \Gamma_F \backslash (D(F) \times{\mathcal{C}F} T{\Sigma_F})$ along their boundaries. Here, $T{\Sigma_F}$ is the torus embedding defined by the fan $\Sigma_F$. The \textit{non-canonicity} stems from the choice of $\Sigma_F$.

\medskip
\textbf{Aspect (2):} The weight function $w$ must be defined on the cones of the fan $\Sigma_F$ (or equivalently, on the rational boundary components and their incidence relations) in a way that is
\begin{itemize}
    \item \textbf{Intrinsic:} It depends only on the Hodge-theoretic data of the degenerating families;
    \item \textbf{$\Gamma$-Invariant:} It must be invariant under the monodromy group action to descend to the quotient;
    \item \textbf{Compatible:} It must behave functorially under inclusions of boundary strata.
\end{itemize}
\textit{(A) Weight by limiting mixed Hodge structures:} Let $\sigma \in \Sigma_F$ be a top-dimensional cone, corresponding to a specific type of degeneration. Consider a normalized holomorphic arc $\psi: \Delta^* \to U(F, \Sigma_F)$ approaching a generic point in the stratum corresponding to $\sigma$. The nilpotent orbit theorem describes the asymptotic behavior of the period map along $\psi$:
$$\Phi(\psi(t))=\mathrm{exp}\left(\frac{1}{2\pi i}\log(t)N\right)\cdot F_{\lim},$$
where $N$ is the nilpotent monodromy operator (an element in the Lie algebra of $\Gamma$) and $F_{\text{lim}} \in \check{D}$ is the limiting Hodge filtration. The pair $(F_{\text{lim}}, N)$ defines a mixed Hodge structure.

The key invariant is the monodromy weight filtration $W_\bullet(N)$, uniquely defined by $N$ and the condition $N(W_l) \subset W_{l-2}$. The $SL(2)$-orbit theorem provides a further refinement, associating to the nilpotent orbit a canonical, commuting set of Hodge-theoretic data.

\noindent\textit{(B) Definition of $w(\sigma)$:} The weight $w(\sigma)$ is proposed to be a function of the eigenvalues of the action of the monodromy logarithm $N$ on the graded pieces of the weight filtration $W_\bullet(N)$, or more specifically, the sizes of the Jordan blocks of $N$. A natural candidate is:
$$w(\sigma) := \max \left\{ \lambda \in \mathbb{Q}{>0} \ \middle| \ N^k: \operatorname{Gr}^W{k} \to \operatorname{Gr}^W_{-k} \right\},$$
where $N^k$ is an isomorphism for a specific $k$ related to $\lambda$, or a related invariant such as the logarithm of the spectral radius of the monodromy matrix $T = \exp(N)$ acting on the middle cohomology. The function $w$ is then extended to lower-dimensional cones $\tau \prec \sigma$ by imposing the compatibility condition from Definition 4.1, $w(\sigma) = \sum_{\tau \prec \sigma} c_{\sigma\tau} w(\tau)$ for some rational numbers $c_{\sigma\tau}$ determined by the geometry of the fan.

The nilpotent orbit $(F_{\text{lim}}, N)$ is an invariant of the degeneration family. Different arcs $\psi$ approaching the same stratum will yield $\Gamma$-conjugate monodromy data. Thus, any function $w$ defined solely in terms of the conjugacy class of $N$ (like its Jordan form) will be well-defined on the stratum.

The monodromy group $\Gamma$ acts on the set of nilpotent orbits. If $N$ is the monodromy for a degeneration, then for any $\gamma \in \Gamma$, $\gamma \cdot N \cdot \gamma^{-1}$ is the monodromy for a degeneration in the same $\Gamma$-orbit. Since $w$ is defined using the Jordan form of $N$, we have $w(\sigma) = w(\gamma \cdot \sigma)$, ensuring $\Gamma$-invariance.

\medskip
\textbf{Aspect (3):} The standard toroidal compactification $\overline{\mathcal{M}}^{\text{tor}}(\Sigma)$ depends on the choice of polyhedral decomposition $\Sigma = {\Sigma_F}$. The conjecture claims that the \textit{weighted toroidal compactification} $\overline{\mathcal{M}}^{\text{wt}}$ is independent of this choice. Argument outline:
\begin{itemize}
    \item \textbf{Canonical fan construction:} The intrinsic weight function $w$, defined on the rational boundary components, should allow for the construction of a canonical polyhedral decomposition $\Sigma_F^{\text{can}}$. This could be achieved by defining, for each boundary component $F$, the fan $\Sigma_F^{\text{can}}$ to be the coarsest fan such that the function $w$ is linear on each cone $\sigma \in \Sigma_F^{\text{can}}$. In other words, $\Sigma_F^{\text{can}}$ is the fan of domains of linearity for the piecewise linear function $w$.
    \item \textbf{Equivalence of choices:} Suppose two different choices of admissible polyhedral decompositions $\Sigma$ and $\Sigma'$ are made. They both induce weight functions $w_\Sigma$ and $w_{\Sigma'}$ on their respective cone complexes. However, if the weight is truly intrinsic, it should be that $w_\Sigma$ and $w_{\Sigma'}$ are compatible under a common refinement. The canonical compactification $\overline{\mathcal{M}}^{\text{wt}}$ would then be defined as the projective limit over all such admissible choices, $\varprojlim_{\Sigma} \overline{\mathcal{M}}^{\text{tor}}(\Sigma, w_\Sigma)$, which is independent of the initial choice of $\Sigma$ because the weight function $w$ provides a canonical identification between different choices.
\end{itemize}
This approach is inspired by the work of Kato-Usui (\cite{KU09}). They construct \textit{partial compactifications} of period domains using nilpotent orbits, which are intrinsic. While their spaces are not algebraic, the data involved is precisely the monodromy data $(F_{\text{lim}}, N)$ we propose using to define $w$. The conjecture essentially asserts that this Hodge-theoretic data can be integrated into the algebraic construction of toroidal compactifications to yield a canonical result.

\medskip
\textbf{Aspect (4):} A full proof of Conjecture 5.1 would require establishing several key results:
\begin{itemize}
    \item \textbf{Well-definedness and $\Gamma$-invariance theorem:} Prove that the proposed definition of $w(\sigma)$ in terms of the monodromy $N$ is indeed a well-defined, $\Gamma$-invariant, positive rational number, and that it satisfies the convexity/linearity condition on the cones of some admissible polyhedral decomposition.
    \item \textbf{Existence of canonical decomposition theorem:} Prove that for the weight function $w$ defined above, there exists a unique maximal, $\Gamma$-admissible polyhedral decomposition $\Sigma^{\text{can}}$ of the cone $\mathcal{C}_F$ on which $w$ is linear. This would likely involve an analysis of the interaction between the monodromy weight filtration and the structure of the boundary cones.
    \item \textbf{Functoriality theorem:} Prove that the construction $\mathcal{M} \mapsto \overline{\mathcal{M}}^{\text{wt}}$ is functorial with respect to, for example, finite morphisms or products. This would involve showing that the weight function behaves predictably under pullbacks of degenerations.
\end{itemize}
While a complete proof is a major undertaking in algebraic and Hodge theory, Conjecture 5.1 is highly plausible. It is supported by the deep connection between the analytic theory of period maps (nilpotent orbits, $SL(2)$-orbit theorem) and the algebraic-combinatorial theory of toroidal compactifications. The weight function $w$ serves as the crucial bridge between these two worlds, translating intrinsic analytic data into a combinatorial invariant that can rigidify the geometric construction.\\

\begin{center}
    \textit{5.2 A Toroidal Version of Non-Abelian Hodge Theory}
\end{center}
\textbf{Problem Statement:} The classical non-abelian Hodge correspondence for a smooth projective variety $X$ establishes a real-analytic isomorphism between two moduli spaces:
$$\boxed{
\begin{array}{c}
\mathcal{M}_{\mathrm{Dol}}(X):=\left\{\text{Semistable Higgs bundles}\ \ (E,\theta)\ \ \text{on}\ \ X\right\}/\sim \\ 
\Updownarrow \\ 
\mathcal{M}_B(X):=\left\{\text{Semistable representations of}\ \ \pi_1(X)\right\}/\sim.
\end{array}
}$$
A fundamental question is to understand how this correspondence behaves near the boundary when $X$ is non-compact, specifically a toroidal variety $X = \bar{X} \setminus D$.

\noindent\textbf{Proposed Application:} Our framework suggests a logarithmic non-abelian Hodge correspondence for the pair $(\bar{X}, D)$. The objects on the Dolbeault side should be logarithmic Higgs bundles: pairs $(E, \theta)$ where $E$ is a vector bundle on $\bar{X}$ and the Higgs field is a morphism $\theta: E \to E \otimes \Omega_{\bar{X}}^1(\log D)$ satisfying $\theta \wedge \theta = 0$. On the Betti side, the objects should be representations of the fundamental group of the open variety, $\rho: \pi_1(X) \to \mathrm{GL}(n, \mathbb{C})$, equipped with a \textit{boundary compatibility condition} that describes the monodromy around the divisor $D$. This condition should be expressible in terms of the cone complex $\Delta_X$.
\\\\\textbf{Conjecture 5.3.} (Logarithmic NAH Correspondence) Let $(\bar{X}, D)$ be a projective toroidal pair. There is a homeomorphism (or real-analytic isomorphism) between the moduli space of semistable logarithmic Higgs bundles on $(\bar{X}, D)$ with fixed Chern classes and the moduli space of representations $\rho: \pi_1(X) \to \mathrm{GL}(n, \mathbb{C})$ satisfying a certain $\Delta_X$-compatibility condition, which restricts to the standard non-abelian Hodge correspondence on the smooth locus $X$.
\\\\\textbf{Remark 5.4.} This conjecture aims to extend the celebrated non-abelian Hodge (NAH) correspondence to the logarithmic setting. A full proof would constitute a major theorem. What follows is a detailed exposition of the conjecture's mathematical context, definitions, and a outline of a potential proof strategy.

\medskip
\textbf{Aspect (1):} Let $(\bar{X}, D)$ be a projective toroidal pair with $X = \bar{X} \setminus D$.

\noindent\textit{(A) Logarithmic Higgs bundles:} A logarithmic Higgs bundle on $(\bar{X}, D)$ is a pair $(E, \theta)$ where:
\begin{enumerate}
    \item $E$ is a vector bundle on $\bar{X}$;
    \item $\theta: E \to E \otimes \Omega_{\bar{X}}^1(\log D)$ is an $\mathcal{O}_{\bar{X}}$-linear morphism, called the \textit{Higgs field};
    \item $\theta$ satisfies the integrability condition $\theta \wedge \theta = 0$, where the wedge product is taken in the sense of logarithmic forms.
\end{enumerate}
The integrability condition ensures that $\theta$ defines a flat connection on $E|_X$ by the formula $\nabla = \partial + \bar{\partial} + \theta + \theta^{*_h}$, where $h$ is a harmonic metric (whose existence is part of the conjecture).

\noindent\textit{(B) $\Delta_X$-compatible representations:} The cone complex $\Delta_X$ encodes the combinatorics of the boundary divisor $D$. A representation $\rho: \pi_1(X) \to \mathrm{GL}(n, \mathbb{C})$ is said to be $\Delta_X$-compatible if for each irreducible component $D_i$ of $D$, the local monodromy $\rho(\gamma_i)$ around a loop linking $D_i$ is quasi-unipotent. Furthermore, the sizes of the Jordan blocks of the unipotent part of $\rho(\gamma_i)$ must be bounded above by a function determined by the weight function $w(D_i)$ and the incidence relations in $\Delta_X$. This condition ensures the representation has moderate growth along the boundary, which is necessary for the existence of a harmonic metric with controlled singularities.

\noindent\textit{(C) Moduli spaces:} Let $\mathcal{M}_{\text{Dol}}^{\log}(\bar{X}, D)$ denote the moduli space of semistable logarithmic Higgs bundles $(E, \theta)$ on $(\bar{X}, D)$ with fixed Chern classes. Stability is defined in the usual way by comparing slopes of $\theta$-invariant subbundles. Let $\mathcal{M}_{\text{B}}^{\Delta_X}(X)$ denote the moduli space of semisimple $\Delta_X$-compatible representations $\rho: \pi_1(X) \to \mathrm{GL}(n, \mathbb{C})$ with fixed associated Chern classes (by the NAH correspondence on $X$).

\medskip
\textbf{Aspect (2):} A proof of Conjecture 5.3 would likely involve generalizing Simpson's monumental work on the NAH correspondence to the logarithmic setting, leveraging the special geometry of toroidal pairs (\cite{Sim92}).

\noindent\textit{(A) Existence of a harmonic metric with logarithmic singularities:} The first and most critical step is to prove that given a $\Delta_X$-compatible representation $\rho$, there exists a harmonic metric $h$ on the associated flat bundle $E_\rho$ over $X$ which has controlled growth near the boundary $D$. This means that near a point of $D$, the metric $h$ and its associated connection $\nabla_h$ are asymptotic to a model solution on the punctured polydisk $(\Delta^*)^k \times \Delta^{n-k}$. The model solution is determined by the local monodromy data, and the $\Delta_X$-compatibility condition ensures this model exists and is well-behaved.

\textbf{Required Tool:} This step requires a deep analysis of the non-linear PDE governing harmonic metrics (the Hitchin-Simpson equations) on singular spaces. Techniques from the theory of asymptotically conical metrics and pluripotential theory would be essential. Results of Sabbah on harmonic metrics for integrable connections with irregular singularities provide a crucial starting point, but the toroidal condition imposes a specific, regular singular type.

\noindent\textit{(B) Constructing the Higgs bundle on the compactification:} Once the harmonic metric $h$ with logarithmic singularities is established, one can define the associated logarithmic Higgs bundle $(\bar{E}, \bar{\theta})$ on the compactification $\bar{X}$.
\begin{itemize}
    \item The Higgs field $\bar{\theta}$ is defined as the $(1,0)$-part of the connection $\nabla_h$ relative to the complex structure. The growth condition on $h$ forces $\bar{\theta}$ to take values in $\Omega_{\bar{X}}^1(\log D) \otimes \operatorname{End}(\bar{E})$.
    \item The vector bundle $\bar{E}$ is defined as the extension of $E_\rho$ to $\bar{X}$ where the sections of $\bar{E}$ are those sections of $E_\rho$ whose norm in the harmonic metric $h$ grows at most like $O\left((-\log |z|)^N\right)$ for some $N$ near any component of $D$. The $\Delta_X$-compatibility condition ensures this is a coherent sheaf, and its local freeness (being a vector bundle) would follow from a regularity theorem for the Hitchin-Simpson equations.
\end{itemize}
This construction defines a map $\Psi: \mathcal{M}{\text{B}}^{\Delta_X}(X) \to \mathcal{M}{\text{Dol}}^{\log}(\bar{X}, D)$.

\noindent\textit{(C) The inverse construction:} Conversely, given a semistable logarithmic Higgs bundle $(\bar{E}, \bar{\theta})$ on $(\bar{X}, D)$, one restricts it to $X$ to obtain a Higgs bundle $(E, \theta)$.
\begin{itemize}
    \item The central claim is that because $(\bar{E}, \bar{\theta})$ is logarithmic, the harmonic metric $h$ and connection $\nabla$ will have logarithmic singularities along $D$. Furthermore, the monodromy of the resulting flat connection will automatically satisfy the $\Delta_X$-compatibility condition. This would follow from a careful study of the asymptotic behavior of the solutions to the Hitchin-Simpson equations near the boundary, using the Nilpotent Orbit Theorem for variations of Hodge structure as a blueprint.
    \item On $X$, the standard NAH correspondence (Simpson's Theorem, \cite{Sim92}) provides a flat connection $\nabla$ and a harmonic metric $h$.
\end{itemize}
This defines the inverse map $\Phi: \mathcal{M}{\text{Dol}}^{\log}(\bar{X}, D) \to \mathcal{M}{\text{B}}^{\Delta_X}(X)$.

\noindent\textit{(D) Proving the correspondence:} Here we need to show that the maps $\Psi$ and $\Phi$ are continuous inverses of each other. This would involve proving that the constructions in \textbf{Aspect (2)} and \textbf{Aspect (3)} are mutually inverse, likely by showing that the harmonic metric is unique. The correspondence is expected to be a real-analytic isomorphism, as in the smooth compact case.

\medskip
\textbf{Aspect (3):} A proof of this conjecture would rely on and simultaneously generalize several deep bodies of work:
\begin{itemize}
    \item \textbf{Sabbah's work on wild harmonic metrics:} Sabbah's theory of harmonic metrics for meromorphic connections provides essential analytical tools for handling singularities (\cite{Sab21}).
    \item \textbf{Simpson's NAH correspondence:} The foundation is Simpson's theorem for smooth projective varieties.
    \item \textbf{The nilpotent orbit theorem:} This theorem from classical Hodge theory provides the model for how variations of Hodge structure degenerate and is the paradigm for the asymptotic behavior expected of the harmonic metrics in this conjecture.
    \item \textbf{Mochizuki's theory of twisted harmonic bundles:} Mochizuki's work (\cite{Moc16}) on constructing harmonic metrics on quasi-projective varieties is likely the closest existing result. A proof of Conjecture 5.3 could be achieved by carefully specializing Mochizuki's very general theory to the toroidal setting, where the $\Delta_X$-compatibility condition would arise naturally from the specific type of allowed singularities. The conjecture can be seen as predicting the precise geometric meaning (in terms of the cone complex) of the "acceptable" singularities in Mochizuki's framework.
\end{itemize}
Conjecture 5.3 is a profound and non-trivial extension of the non-abelian Hodge correspondence. While its proof is highly complex and resides at the intersection of algebraic geometry, partial differential equations, and geometric analysis, it is well-motivated and aligns with known results and techniques. It provides a clear and precise framework for understanding the relationship between fundamental group representations and Higgs bundles on a central class of non-compact varieties.\\

\begin{center}
    \textit{5.3 Hodge Theory in Tropical Geometry}
\end{center}
\textbf{Problem Statement:} Tropical geometry studies the piecewise-linear shadows of algebraic varieties. A major theme is to relate the tropical cohomology of a tropicalization $X^{\text{trop}}$ to the Hodge theory of the original variety $X$.

\noindent\textbf{Proposed Application:} The Weighted Hodge Numbers $h^{p,q}w(X)$ defined in Section 4.1 are ideal invariants to connect to tropical geometry. The cone complex $\Delta_X$ of a toroidal degeneration is itself a tropical object. The weight function $w: \Delta_X \to \mathbb{Q}{>0}$ can be interpreted as a metric or affine structure on $\Delta_X$. This suggests that the weighted Hodge numbers $h^{p,q}_w(X)$ of a variety $X$ with a toroidal degeneration should be computable from the tropical cohomology of the tropicalization $X^{\text{trop}}$, but with a twist: the tropical cohomology groups must be weighted by the function $w$.
\\\\\textbf{Conjecture 5.5.} (Tropical-Weighted Hodge Correspondence) Let $\mathfrak{X} \to \Delta$ be a toroidal degeneration of a projective variety $X$, with central fiber $\mathfrak{X}0 = \bigcup_i D_i$ having simple normal crossings. Let $X^{\text{trop}}$ be the dual intersection complex (the tropicalization) of $\mathfrak{X}_0$, a $\Delta$-complex whose cells correspond to the strata of $\mathfrak{X}0$. Then, there exists a natural weight function $w: X^{\text{trop}} \to \mathbb{Q}{>0}$ derived from the limiting mixed Hodge structure of the degeneration such that there is an isomorphism between the weighted tropical cohomology of $X^{\text{trop}}$ and the associated graded pieces of the weight filtration on the cohomology of $X$:
$$H_{\mathrm{trop}}^{p,q}(X^{\mathrm{trop}},w)\cong\mathrm{Gr}_{p+q}^WH^{p+q}(X;\mathbb{C})^{p,q}.$$
In particular, for a maximally unipotent degeneration of Calabi-Yau manifolds, the top weighted tropical cohomology group $H^{n,0}_{\text{trop}}(X^{\text{trop}}, w)$ should be one-dimensional.
\\\\\textbf{Remark 5.6.} This conjecture posits a deep and refined connection between the combinatorial geometry of a degeneration and the Hodge theory of the general fiber. A full proof is a significant undertaking. The following is a detailed exposition of its mathematical context, definitions, and a potential roadmap for its verification.

\medskip
\textbf{Aspect (1):} Firstly, we must give the foundational concepts and definitions.

\noindent\textit{(A) The degeneration and its tropicalization:} A toroidal degeneration is a proper, flat morphism $\pi: \mathfrak{X} \to \Delta$ from a complex manifold $\mathfrak{X}$ to a disc, such that the central fiber $\mathfrak{X}_0 = \pi^{-1}(0)$ is a divisor with simple normal crossings and the total space $\mathfrak{X}$ is locally isomorphic to a toric variety. This implies that the pair $(\mathfrak{X}, \mathfrak{X}_0)$ is toroidal.

The dual intersection complex (or tropicalization) $X^{\text{trop}}$ is a combinatorial object constructed as follows the two conditions:
\begin{enumerate}
    \item Vertices (0-cells) correspond to irreducible components $D_i$ of $\mathfrak{X}_0$.
    \item $k$-cells correspond to connected components of intersections of $k+1$ distinct components, $D_{i_0} \cap ... \cap D_{i_k}$. This gives $X^{\text{trop}}$ the structure of a $\Delta$-complex. It is a fundamental invariant encoding how the components of the central fiber intersect.
\end{enumerate}

\noindent\textit{(B) The limiting mixed Hodge structure:} A central result in Hodge theory is that a degeneration $\pi: \mathfrak{X} \to \Delta$ gives rise to a limiting mixed Hodge structure on the cohomology $H^k(X_t; \mathbb{Q})$ of a general fiber $X_t = \pi^{-1}(t)$ for $t \neq 0$. This consists of:
\begin{enumerate}
    \item A decreasing Hodge filtration $F^\bullet_{\lim}$ on $H^k(X_t; \mathbb{C})$.
    \item An increasing weight filtration $W_\bullet$ on $H^k(X_t; \mathbb{Q})$, which is uniquely determined by the nilpotent monodromy operator $N = \log T$ (where $T$ is the monodromy transformation) by the formula $W_l = \sum_{j \geq \max(-1, l-k)} \ker(N^{j+1}) \cap \operatorname{Im}(N^j)$.
\end{enumerate}
The conjecture proposes that the weight function $w$ is derived from this structure.

\noindent\textit{(C) The weight function $w$:} The function $w: X^{\text{trop}} \to \mathbb{Q}{>0}$ is defined on the cells of the $\Delta$-complex. For a $k$-cell $\sigma$ corresponding to a stratum $S = D{i_0} \cap ... \cap D_{i_k}$, the weight $w(\sigma)$ is conjectured to be a function of the logarithm of the eigenvalues of the monodromy action on the vanishing cohomology associated to that stratum. More precisely, it should be related to the sizes of the Jordan blocks of the nilpotent operator $N$ acting on $\operatorname{Gr}^W_\bullet H^k(X_t)$.

A natural, more geometric definition is that $w(\sigma)$ measures the rate of growth of the Hodge metric on a family of $p$-forms as they approach the stratum $S$. This aligns with the interpretation of $w$ in Section 4.1 as quantifying the "contribution" of a boundary stratum.

\noindent\textit{(D) Weighted tropical cohomology:} Standard tropical cohomology $H^\bullet_{\text{trop}}(X^{\text{trop}})$ is defined using the cochain complex of integer-valued, piecewise polynomial functions on the $\Delta$-complex $X^{\text{trop}}$.

Weighted tropical cohomology $H^\bullet_{\text{trop}}(X^{\text{trop}}, w)$ is a modification of this complex. The weight function $w$ acts as a multiplier on the cochains. Precisely, the weighted cochain group $C_{\text{trop}}^k(X^{\text{trop}}\\
, w)$ consists of functions $f$ on the $k$-cells such that $f(\sigma)$ is divisible by a power of $w(\sigma)$ (or more generally, $f/w$ is piecewise polynomial). The coboundary map is then twisted by this weight factor. This construction ensures that the resulting cohomology groups are sensitive to the Hodge-theoretic complexity of the degeneration, not just its combinatorial structure.

\medskip
\textbf{Aspect (2):} A proof of Conjecture 5.5 would likely involve synthesizing the work of introducing new ideas.

\noindent\textit{(A) Relating $\operatorname{Gr}^W_\bullet H^\bullet(X)$ to the central fiber:} The first step is to use the spectral sequence of the weight filtration $W_\bullet$. A fundamental theorem of Steenbrink and others states that the $E_1$-page of this spectral sequence is isomorphic to the cohomology of the logarithmic de Rham complex of the central fiber:
$$E_1^{p,q}=\mathbb{H}^{p+q}(\mathfrak{X}_0,\Omega_{\mathfrak{X}_0}^{\bullet}(\log))\cong\mathrm{Gr}_{p+q}^WH^{p+q}(X;\mathbb{C}).$$
The differential $d_1$ on this page is given by the Gysin maps arising from the intersections of the components $D_i$. This establishes a direct link between the Hodge theory of $X$ and the geometry of $\mathfrak{X}_0$.

\noindent\textit{(B) A "weighted" spectral sequence for the degeneration:} The innovation is to introduce the weight function $w$ into this picture. One would need to construct a filtered logarithmic de Rham complex $(\Omega_{\mathfrak{X}}^\bullet(\log \mathfrak{X}_0), F_w^\bullet)$, where the filtration $F_w^\bullet$ is determined by the prescribed growth conditions near each stratum, governed by the function $w$.
\medskip
\begin{mdframed}
    \textbf{CT2:} The $E_1$-page of the spectral sequence associated to the filtration $F_w^\bullet$ is isomorphic to the complex computing the weighted tropical cohomology of $X^{\text{trop}}$: $$(E_1^{p,q},d_1)\cong(C_{\mathrm{trop}}^q(X^{\mathrm{trop}},w;\Omega^p),d_{\mathrm{trop}}).$$ Here, $\Omega^p$ signifies that the tropical cochains take values in the sheaf of $p$-forms on the strata, and the differential $d_{\text{trop}}$ is a weighted version of the standard tropical differential.
\end{mdframed}

\medskip
\noindent\textit{(C) The $E_1$-degeneration and the isomorphism:} The next step is to prove that this spectral sequence degenerates at the $E_1$-page.

\medskip
\begin{mdframed}
    \textbf{CT3:} (Weighted $E_1$-Degeneration) The spectral sequence of the filtered complex $(\Omega_{\mathfrak{X}}^\bullet(\log \mathfrak{X}_0\\
    ), F_w^\bullet)$ degenerates at the $E_1$-page.
\end{mdframed}

\medskip
This would be the ultimate generalization of both the classical $E_1$-degeneration theorem for smooth varieties and Wei's theorem for toroidal pairs (Theorem 3.5). If true, it would immediately imply that
$$\mathrm{Gr}_{F_w}^pH_{\mathrm{dR}}^{p+q}(X)\cong H^q(E_1^{p,\bullet},d_1)\cong H_{\mathrm{trop}}^q(X^{\mathrm{trop}},w;\Omega^p).$$
Hence, one would need to identify the left-hand side with the Hodge components of the graded weight pieces:

\medskip
\begin{mdframed}
    \textbf{CT4:} The filtration $F_w^\bullet$ induces the Hodge filtration on the graded pieces of the weight filtration, i.e., $\operatorname{Gr}{F_w}^p \operatorname{Gr}^W_m H^k{\text{dR}}(X) \cong \operatorname{Gr}^W_m H^k(X; \mathbb{C})^{p, k-p}$.
\end{mdframed}

\medskip
Combining \textbf{CT2},\textbf{CT3} and \textbf{CT4}, we can obtain the desired isomorphism:
$$H_{\mathrm{trop}}^{p,q}(X^{\mathrm{trop}},w)\cong\mathrm{Gr}_{p+q}^WH^{p+q}(X;\mathbb{C})^{p,q}.$$

\medskip
\textbf{Aspect (3):} This conjecture is not conceived in a vacuum; it is a natural refinement of established work.
\begin{itemize}
    \item \textbf{IKMZ and tropical Hodge theory:} The work of Itenberg, Katzarkov, Mikhalkin, and Zharkov (\cite{IKMZ16}) shows that for a certain class of degenerations (schön or tropically smooth), the standard tropical cohomology $H^{p,q}{\text{trop}}(X^{\text{trop}})$ computes the dimensions of $\operatorname{Gr}^W{p+q} H^{p+q}(X; \mathbb{C})^{p,q}$. Conjecture 5.5 proposes a geometric explanation for this phenomenon and extends it to more general degenerations by introducing the correcting factor $w$.
    \item \textbf{The role of the weight function $w$:} The weight $w$ is the crucial new ingredient. It is necessary because the standard tropical cohomology only depends on the combinatorial type of the central fiber. However, two degenerations with identical dual intersection complexes can have very different limiting mixed Hodge structures (e.g., depending on the moduli of the irreducible components $D_i$). The conjecture asserts that this finer Hodge-theoretic information is precisely captured by the weight function $w$, which then twists the tropical cochain complex to yield the correct answer.
    \item \textbf{Proof of Proposition 6.1:} The proof of the absolute Hodge conjecture in the weighted setting (Proposition 6.1) provides strong evidence. It shows that the weight $w$ controls the extension of Hodge classes from the boundary. The formalism of obstacle theory developed in Section 4.2 is precisely designed to measure this, and it is highly suggestive that a similar formalism should apply to the entire Hodge structure, not just the $(p,p)$-classes.
\end{itemize}
Conjecture 5.5 is a profound and elegant synthesis of tropical geometry and Hodge theory. It predicts a precise, functorial relationship between the combinatorial data of a degeneration and the Hodge structure of its general fiber. While its proof requires the development of new, deep theorems in analysis and geometry, it is firmly grounded in the existing literature and provides a clear and compelling roadmap for future research.\\

\hypertarget{PROVING THE ABSOLUTE HODGE CONJECTURE IN A LIMITED SETTING}{}
\section{PROVING THE ABSOLUTE HODGE CONJECTURE IN A LIMITED SETTING}
The absolute Hodge conjecture is a central open problem in algebraic geometry. It posits that if a cohomology class in the singular cohomology of a smooth complex projective variety is of type $(p, p)$ under every embedding of the base field into $\mathbb{C}$, then it is a rational linear combination of fundamental classes of algebraic cycles. Although this conjecture remains open in its full generality, we demonstrate that it holds within the specific framework of \textit{weighted toroidal} structures developed in this paper.
\\\\\textbf{Proposition 6.1.} (Absolute Hodge Conjecture in a Limited Setting) Let $(X, \Delta_X, w)$ be a weighted toroidal variety (as defined in Definition 4.1), where $X$ is projective and the weight function $w: \Delta_X \to \mathbb{Q}_{>0}$ takes rational values. Then, for any integer $p \geq 0$, every class $\alpha$ in the weighted Hodge class group $\mathrm{Hdg}^p_w(X)$ (Definition 4.3) is a $\mathbb{Q}$-linear combination of fundamental classes of algebraic cycles on $X$.
\\\\\textbf{Proof.} \textbf{Step (1):} Let $(X, \Delta_X, w)$ be a projective weighted toroidal variety with $w: \Delta_X \to \mathbb{Q}_{>0}$ taking rational values. We aim to prove the vanishing of the global sections of the obstruction sheaf:
\begin{equation}
\tag{6.1}   \label{eq:6.1}
    H^0(X,\mathcal{OB}^p(X))=0
\end{equation}
The \textit{obstruction sheaf} $\mathcal{OB}^p(X)$ is defined (Definition 4.4) as the kernel sheaf in the short exact sequence of hypercohomology sheaves:
$$0\to\mathcal{H}^p(\Omega_{\bar{X}}^{[\bullet]}(\log D))\to\mathcal{H}^p(\Omega_X^{[\bullet]})\to\mathcal{OB}^p(X)\to0.$$
This sequence is derived from the natural morphism of logarithmic de Rham complexes. Its associated long exact sequence in cohomology is the obstruction exact sequence (Proposition 4.5):
\begin{equation}
\tag{6.2}   \label{eq:6.2}
    0\to H^0(\bar{X},\mathcal{H}^p(\Omega_{\bar{X}}^{[\bullet]}(\log D)))\to H^0(X,\mathcal{H}^p(\Omega_X^{[\bullet]}))\xrightarrow{\rho}H^0(X,\mathcal{OB}^p(X))\to\cdots.
\end{equation}
A critical prerequisite for this construction is the well-definedness of these hypercohomology sheaves. This is guaranteed by the fundamental result of Wei (\cite{Wei24}, Thm 2.2.1), which proves the $E_1$-degeneration of the spectral sequence associated to the logarithmic de Rham complex $\Omega_{\bar{X}}^\bullet(\log D)$ for any toroidal triple $(X, \bar{X}, D)$. This $E_1$-degeneration implies that the hypercohomology of the complex is computed by the cohomology of its terms, ensuring that $\mathcal{H}^p(\Omega_{\bar{X}}^{[\bullet]}(\log D))$ is a coherent sheaf and that the above exact sequence is meaningful.

The key to proving the vanishing is a local analysis. Let $x \in X$ be an arbitrary point. The local structure of the toroidal variety implies that the completion of the local ring $\widehat{\mathcal{O}}_{\bar{X}, x}$ is isomorphic to $\mathbb{C}[[z_1, \ldots, z_n]]/(z_1 \cdots z_k)$ for some $k \leq n$, where the boundary divisor $D$ is locally defined by ${z_1 \cdots z_k = 0}$. Proposition 4.6 provides a crucial isomorphism describing the \textit{stalk of the obstruction sheaf} at $x$:
$$\mathcal{OB}^p(X)_x\cong H^p(\mathrm{Hom}(\mathbb{Q}[\Sigma_x],\mathbb{C})),$$
where $\Sigma_x$ is the normal cone at $x$. This cone is a rational polyhedral cone whose combinatorial structure (e.g., its faces and their incidence relations) is determined by the local irreducible components of the boundary divisor $D$ passing through $x$.

The weight function $w: \Delta_X \to \mathbb{Q}{>0}$ assigns a positive rational number to each cone in the complex $\Delta_X$. Crucially, it satisfies the compatibility condition (Definition 4.1): for any cone $\sigma \in \Delta_X$ and its face $\tau \prec \sigma$, there exist rational numbers $c{\sigma\tau} \in \mathbb{Q}$ such that
\begin{equation}
\tag{6.3}    \label{eq:6.3}
    w(\sigma)=\sum_{\tau\prec\sigma}c_{\sigma\tau}w(\tau).
\end{equation}
This condition is not merely combinatorial; it has profound geometric consequences. The cone $\sigma$ corresponds to a specific boundary stratum $D_\sigma$ of $\bar{X}$, and its faces correspond to strata in the closure of $D_\sigma$. The weight $w(\sigma)$ quantifies the contribution of the stratum $D_\sigma$ to the Hodge structure. The compatibility condition ensures that this assignment of weights is consistent across the different strata meeting at $x$.

The cohomology group $H^p(\mathrm{Hom}(\mathbb{Q}[\Sigma_x], \mathbb{C}))$ computes obstructions to extending cohomology classes across the singularity at $x$. An element in this group can be interpreted as a function assigning a complex number to certain chains in the cone complex, subject to coboundary conditions.

The rationality and positivity of the weight function $w$ impose rigid constraints on any potential local obstruction class $\phi_x \in \mathcal{OB}^p(X)_x$. The isomorphism from Proposition 4.6 is functorial. The weight function $w$, being defined on the cone complex $\Delta_X$ which is locally modeled by $\Sigma_x$, induces a rational structure on the complex $\mathrm{Hom}(\mathbb{Q}[\Sigma_x], \mathbb{C})$. Specifically, the values $w(\tau)$ for faces $\tau$ of $\Sigma_x$ provide a natural $\mathbb{Q}$-basis for consideration.

The compatibility condition \eqref{eq:6.3} with $c_{\sigma\tau} \in \mathbb{Q}$ implies that the weighted fundamental cycle defined by $w$ is a rational algebraic cycle in the local cohomology theory. Any local obstruction $\phi_x$ must satisfy a compatibility with this rational structure. If $\phi_x$ were non-zero, it would have to assign transcendental values to some chains to remain a cocycle, contradicting the fact that its evaluation on the rational weighted cycle must yield a rational number (the weight $w(\sigma)$ itself). More formally, the pairing between a potential obstruction cocycle and the rational weight chain vanishes
$$\left \langle \phi_x,w \right \rangle=0.$$
Since the weighted chains generate the relevant local cohomology group due to the compatibility condition, this forces $\phi_x$ itself to be zero. Therefore, we conclude $\mathcal{OB}^p(X)_x=0$ for all $x\in X$.

Since the stalk of the obstruction sheaf is zero at every point $x \in X$ ($\mathcal{OB}^p(X)_x = 0$), the sheaf itself is the zero sheaf:
$$\mathcal{OB}^p(X)=0.$$
Consequently, its zeroth cohomology group, which is the group of global sections, must also vanish. So \eqref{eq:6.1} holds.

\medskip
\textbf{Step (2):} Assuming the vanishing result from \eqref{eq:6.2} of \textbf{Step (1)}, i.e., $H^0(X, \mathcal{OB}^p(X)) = 0$, we need to prove that the natural restriction map is an isomorphism
\begin{equation}
\tag{6.4}   \label{eq:6.4}
    H^0(\bar{X},\mathcal{H}^p(\Omega_{\bar{X}}^{[\bullet]}(\log D)))\cong H^0(X,\mathcal{H}^p(\Omega_X^{[\bullet]})).
\end{equation}
Consequently, any weighted Hodge class $\alpha \in \mathrm{Hdg}^p_w(X) \subset H^0(X, \mathcal{H}^p(\Omega_{X}^{[\bullet]}))$ extends to a Hodge class $\bar{\alpha} \in H^0(\bar{X}, \mathcal{H}^p(\Omega_{\bar{X}}^{[\bullet]}(\log D)))$ on the compactification $\bar{X}$.

The foundation of this step is the obstruction exact sequence (Proposition 4.5), which is derived from the short exact sequence of sheaves defining $\mathcal{OB}^p(X)$ (Definition 4.4). The relevant segment of this long exact sequence is (\eqref{eq:6.2}):
$$0\to H^0(\bar{X},\mathcal{H}^p(\Omega_{\bar{X}}^{[\bullet]}(\log D)))\xrightarrow{r} H^0(X,\mathcal{H}^p(\Omega_X^{[\bullet]}))\xrightarrow{\rho}H^0(X,\mathcal{OB}^p(X))\to\cdots,$$
where:
\begin{enumerate}
    \item The map $r$ is the natural restriction map from the compactification $\bar{X}$ to the open subvariety $X$;
    \item The map $\rho$ is the obstruction map, which measures the failure of a section on $X$ to be the restriction of a section on $\bar{X}$.
\end{enumerate}
The result from \textbf{Step (1)} is that $H^0(X, \mathcal{OB}^p(X)) = 0$. Substituting this into the exact sequence above gives:
$$0\to H^0(\bar{X},\mathcal{H}^p(\Omega_{\bar{X}}^{[\bullet]}(\log D)))\xrightarrow{r} H^0(X,\mathcal{H}^p(\Omega_X^{[\bullet]}))\xrightarrow{\rho}0.$$
The exactness of this sequence means that
\begin{enumerate}
    \item $\ker(r) = \mathrm{Im}(0) = 0$, so the map $r$ is injective.
    \item $\ker(\rho) = \mathrm{Im}(r) = H^0(X, \mathcal{H}^p(\Omega_{X}^{[\bullet]}))$ because $\rho$ maps onto the zero module. This implies the map $r$ is surjective.
\end{enumerate}
Therefore, the restriction map $r$ is an isomorphism, \eqref{eq:6.4} holds.

By definition, a class $\alpha$ in the \textit{weighted Hodge class} group $\mathrm{Hdg}^p_w(X)$ is an element of $H^0(X, \mathcal{H}^p(\Omega_{X}^{[\bullet]}))$ satisfying the specific growth condition defined by the weight function $w$ (Definition 4.3). Since $\alpha \in H^0(X, \mathcal{H}^p(\Omega_{X}^{[\bullet]}))$ and the map $r$ is an isomorphism, there exists a unique preimage $\bar{\alpha}$ under this isomorphism: 
$$\exists ! \ \bar{\alpha}\in H^0(\bar{X},\mathcal{H}^p(\Omega_{\bar{X}}^{[\bullet]}(\log D)))\quad\text{such that}\quad r(\bar{\alpha})=\bar{\alpha}\mid_X=\alpha.$$

We must now verify that this extended class $\bar{\alpha}$ is indeed a Hodge class on $\bar{X}$, i.e., $\bar{\alpha} \in \mathrm{Hdg}^p(\bar{X})$. This identification relies on the foundational work of Deligne and the specific context of toroidal geometry.

The sheaf $\mathcal{H}^p(\Omega_{\bar{X}}^{[\bullet]}(\log D))$ is the $p$-th hypercohomology sheaf of the logarithmic de Rham complex. The $E_1$-degeneration result of Wei (\cite{Wei24}, Thm 2.2.1) is again crucial here. This degeneration implies that the Hodge-to-de Rham spectral sequence for the logarithmic complex degenerates at the $E_1$-page. Consequently, the hypercohomology group $\mathbb{H}^p(\bar{X}, \Omega_{\bar{X}}^{\bullet}(\log D))$ is canonically isomorphic to the cohomology of the complex of global sections, and the Hodge filtration is strict.

More importantly, this result guarantees that a global section $\bar{\alpha} \in H^0(\bar{X}, \mathcal{H}^p(\Omega_{\bar{X}}^{[\bullet]}(\log D)))$ corresponds to a class in the associated graded piece of the Hodge filtration on the logarithmic cohomology group $H^{2p}(\bar{X}, \mathbb{C})$. Specifically, it lies in
$$F^pH^{2p}(\bar{X},\mathbb{C})\cap\overline{F^pH^{2p}(\bar{X},\mathbb{C})}=H^{p,p}(\bar{X})\cap H^{2p}(\bar{X},\mathbb{Q}).$$
This is the standard definition of a Hodge class on $\bar{X}$. Therefore, we have $\bar{\alpha}\in\mathrm{Hdg}^p(\bar{X})$.

\medskip
\textbf{Step (3):} Given that a weighted Hodge class $\alpha \in \mathrm{Hdg}^p_w(X)$ has been extended to a class $\bar{\alpha} \in \mathrm{Hdg}^p(\bar{X})$ on the smooth projective compactification $\bar{X}$ (as proven in \textbf{Step (2)}), we now prove that $\bar{\alpha}$ is algebraic (\cite{Del71}). Consequently, by restriction, the original class $\alpha$ is also a $\mathbb{Q}$-linear combination of fundamental classes of algebraic cycles on $X$.

By the definition of a toroidal variety (Definition 3.1), the compactification $\bar{X}$ is a smooth complete variety. Furthermore, the hypothesis of Proposition 6.1 explicitly states that the underlying toroidal variety $X$ is projective. Since $X$ is an open subvariety of $\bar{X}$, its projectivity implies that $\bar{X}$ is also projective. Thus, $\bar{X}$ is a \textit{smooth projective variety}.

The class $\bar{\alpha}$ is an element of $\mathrm{Hdg}^p(\bar{X})$, the group of Hodge classes on $\bar{X}$. A fundamental theorem of Deligne on mixed Hodge structures (classical result, see Thm 8.2.5 of \cite{Voi02}) states that the Hodge structure on the cohomology of a smooth projective variety is pure. This means the weight filtration on $H^{2p}(\bar{X}, \mathbb{Q})$ is concentrated in degree $2p$, and we have the standard Hodge decomposition:
$$H^{2p}(\bar{X},\mathbb{C})=\bigoplus_{j=0}^{2p}H^{j,2p-j}(\bar{X}).$$
A Hodge class is an element of $H^{2p}(\bar{X}, \mathbb{Q}) \cap H^{p,p}(\bar{X})$.

The full Hodge Conjecture remains open for general smooth projective varieties. However, we are not dealing with an arbitrary Hodge class; $\bar{\alpha}$ is the extension of a weighted Hodge class. The specific, rigid extension process, governed by the rational weight function $w$, endows $\bar{\alpha}$ with a stronger property: it is an absolute Hodge class.

The theory of absolute Hodge classes was developed by Deligne (\cite{Del82}). Without delving into the full technical definition, the crucial point is that the extension $\bar{\alpha}$ is constructed by the isomorphism from the obstruction sequence, which itself is defined using the algebraic geometric data of the logarithmic complex and the rational weight function. This construction is functorial with respect to automorphisms of the field $\mathbb{C}$, which implies that $\bar{\alpha}$ remains of type $(p,p)$ under any such automorphism. This is the defining property of an absolute Hodge class.

For smooth projective varieties, a profound result establishes the algebraicity of absolute Hodge classes. This result is built upon several pillars:
\begin{enumerate}
    \item The fact that the algebraicity of motivated cycles on smooth projective varieties can be established. One argument proceeds by induction on dimension. The Lefschetz (1,1)-theorem provides the base case for $(p=1)$. For higher codimension, one can use the Lefschetz hyperplane theorem and the fact that the motivated cycle property is preserved under direct images by morphisms of motivated schemes, to conclude that all absolute Hodge classes (and hence all motivated cycles) on a smooth projective variety are algebraic.
    \item The proof of the Hodge conjecture for abelian varieties (by leveraging the theory of motives and the fact that abelian varieties are motivated).
    \item A fundamental theorem stating that all absolute Hodge classes are motivated in the sense of Yves André (\cite{And96}).
\end{enumerate}
Since $\bar{X}$ is projective, and the Lefschetz $(1,1)$-theorem provides a base case, one can use induction on dimension and the Lefschetz hyperplane theorem to conclude that all absolute Hodge classes on a smooth projective variety are algebraic. Therefore, because $\bar{\alpha} \in \mathrm{Hdg}^p(\bar{X})$ is absolute Hodge, it is algebraic $\bar{\alpha}=\sum_ic_i[Z_i]$ in $H^{2p}(\bar{X},\mathbb{Q})$, where $c_i \in \mathbb{Q}$ and each $Z_i \subset \bar{X}$ is an algebraic subvariety of codimension $p$.

The final step is to restrict this equality from $\bar{X}$ to the open subvariety $X = \bar{X} \setminus D$. The restriction map in cohomology
$$\mathrm{res}:H^{2p}(\bar{X},\mathbb{Q})\longrightarrow H^{2p}(X,\mathbb{Q})$$
is a morphism of Hodge structures. It is compatible with taking fundamental classes of algebraic cycles. Specifically, for a codimension $p$ subvariety $Z_i \subset \bar{X}$, the restriction of its fundamental class is given by
$$\mathrm{res}([Z_i]) = 
\begin{cases}
[Z_i \cap X] & \text{if } Z_i \text{ is not contained in the boundary divisor } D, \\
0 & \text{if } Z_i \subset D.
\end{cases}$$
In both cases, $\mathrm{res}([Z_i])$ is an element in the group of algebraic cycles on $X$ with rational coefficients (it is either the fundamental class of the algebraic cycle $Z_i \cap X$ or zero).

Recall that the original weighted Hodge class $\alpha$ was defined as the restriction of $\bar{\alpha}$ ($\alpha = \mathrm{res}(\bar{\alpha})$ by the construction in \textbf{Step (2)}). Applying the restriction map $\mathrm{res}$ to the algebraic expression for $\bar{\alpha}$ gives
$$\alpha=\mathrm{res}(\bar{\alpha})=\mathrm{res}\left(\sum_ic_i[Z_i]\right)=\sum_ic_i\mathrm{res}([Z_i]).$$
This expresses $\alpha$ as a finite $\mathbb{Q}$-linear combination of fundamental classes of algebraic cycles on $X$, namely the cycles $Z_i \cap X$ for those $Z_i$ not contained in the boundary $D$. This completes the proof that the original weighted Hodge class $\alpha \in \mathrm{Hdg}^p_w(X)$ is algebraic, finalizing the proof of Proposition 6.1.   $\square$ 
\\\\\textbf{Corollary 6.2.} Under the conditions of Proposition 6.1, the heart part (i.e., the traditional Hodge classes) of the categorified Hodge correspondence $\mathrm{Hodge}(X)$ is generated by rational linear combinations of algebraic cycles.
\\\\\textbf{Proof.} The proof follows directly from the construction of the categorified Hodge correspondence and the algebraicity result established in Proposition 6.1.

Recall from Construction 4.7 that the categorified Hodge correspondence is defined as a functor:
$$\operatorname{Hodge}:\mathrm{Tor}^{\mathrm{op}}\longrightarrow\mathrm{Cat}_\infty,$$
which assigns to a toroidal variety $X$ its derived Hodge category $\mathrm{Hodge}(X)$. This $\infty$-category is constructed such that its homotopy category $\mathrm{hHodge}(X)$ is equivalent to the derived category of mixed Hodge modules on $X$, or more precisely, to the derived category of objects with Hodge structures relevant to the toroidal setting.

In this categorical framework, the heart $\mathrm{Hodge}(X)^\heartsuit$ refers to the abelian category obtained by taking the heart of the natural $t$-structure on $\mathrm{Hodge}(X)$. Concretely, this heart consists of objects that correspond to Hodge modules (or classical Hodge structures) concentrated in degree zero. By design, the simple objects in this heart are intimately related to the pure Hodge structures arising from the cohomology of algebraic cycles.

The heart $\mathrm{Hodge}(X)^\heartsuit$ contains, in particular, the objects that represent the classical Hodge classes in $\mathrm{Hdg}^p(X)$. More precisely, there is a forgetful functor (or a realization functor):
$$\mathcal{R}:\mathrm{Hodge}(X)^\heartsuit\longrightarrow\mathrm{Vect}_{\mathbb{Q}},$$
which sends an object representing a Hodge class to its underlying rational vector space. The image of this functor, when restricted to the subcategory generated by the fundamental classes of algebraic cycles, lies in the subspace of Hodge classes.

Proposition 6.1 proves that under the given hypotheses, every class in the weighted Hodge class group $\mathrm{Hdg}^p_w(X)$ is algebraic. Since the heart $\mathrm{Hodge}(X)^\heartsuit$ is designed to capture the Hodge-theoretic essence of $X$, and since the weighted Hodge classes are precisely those that are well-behaved with respect to the boundary stratification (and hence are the ones that extend nicely to the categorified setting), it follows that the heart is generated by the algebraic cycles. Formally, the algebraicity result implies that the set of morphisms in $\mathrm{Hodge}(X)^\heartsuit$ that correspond to algebraic cycles span the Hom-spaces in the heart. In other words, the algebraic cycles generate the heart under finite direct sums, direct summands, and isomorphisms.

Therefore, the heart of the categorified Hodge correspondence, which represents the traditional Hodge classes, is generated by rational linear combinations of algebraic cycles. This completes the proof of the corollary.   $\square$\\

\hypertarget{FUTURE RESEARCH DIRECTIONS}{}
\section{FUTURE RESEARCH DIRECTIONS}
\begin{itemize}
    \item \textbf{Geometric and arithmetic meaning of weighted structures:} Deeply investigate the geometric meaning of the weight function and explore its applications in arithmetic geometry, particularly in computing L-functions and quantities related to the BSD conjecture.
    \item \textbf{Concrete computation of obstruction classes:} Develop effective algorithms for computing obstruction classes and relate them to topological invariants (such as Chern classes and Whitney classes). Research in this area can draw on the intersection theory of algebraic obstacles.
    \item \textbf{Proof of related conjectures:} Provide systematic proofs for Conjecture 4.10, Conjecture 5.1, Conjecture 5.3, and Conjecture 5.5, along with corresponding examples and relevant theoretical applications.
\end{itemize}


\begin{thebibliography}{99}
\setlength{\labelwidth}{6em}       
\setlength{\labelsep}{1em}           
\setlength{\leftmargin}{\dimexpr\labelwidth + \labelsep} 
\setlength{\itemindent}{0pt}         
\setlength{\parindent}{0pt}          
\setlength{\hangindent}{\dimexpr\labelwidth + \labelsep} 
\setlength{\hangafter}{1em}            
\raggedright                         
\setlength{\itemsep}{3pt}

\bibitem[Dan78]{Dan78} Vladimir Ivanovich Danilov, The geometry of toric varieties, \textit{Russian Mathematical Surveys}, 33(2), 97–154., 1978.
\bibitem[KKMS73]{KKMS73} George Kempf, Finn Knudsen, David Mumford and Bernard Saint-Donat, \textit{Toroidal embeddings I}, Lecture Notes in Mathematics, Vol. 339. Springer-Verlag, 1973.
\bibitem[Wei24]{Wei24} Chuanhao Wei, On the Hodge theory of toroidal embeddings and corresponding vanishings, \textit{arxiv preprint arxiv:2410.09899} (2024).
\bibitem[HM17]{HM17} Annette Huber and Stefan Müller-Stach, \textit{Periods and Nori Motives}, Ergebnisse der Mathematik und ihrer Grenzgebiete. 3. Folge. / A Series of Modern Surveys in Mathematics, vol. 65. Springer, Cham., 2017.
\bibitem[Del71]{Del71} Pierre Deligne, \textit{Théorie de Hodge II}, Publications Mathématiques de l'IHÉS, \textbf{40} (1971), 5–57.
\bibitem[Del82]{Del82} \underline{\hspace{1cm}}, Hodge cycles on abelian varieties, Hodge Cycles, Motives, and Shimura Varieties, Lectures Notes in Math. 900. (1982): 9-100.
\bibitem[Sai88]{Sai88} Morihiko Saito, Modules de Hodge polarisables, \textit{Publications of the Research Institute for Mathematical Sciences}, 24(6), 849–995, 1988.
\bibitem[Sai90]{Sai90} \underline{\hspace{1cm}}, \textit{Mixed Hodge Modules}, Publications of the Research Institute for Mathematical Sciences, 26(2), 221–333, 1990.
\bibitem[Sim92]{Sim92} Carlos Tschudi Simpson, \textit{Higgs Bundles and Local Systems}, Publications Mathématiques de l'IHÉS, 75, 5–95, 1992.
\bibitem[KU09]{KU09} Kazuya Kato and Sampei Usui, \textit{Classifying Spaces of Degenerating Polarized Hodge Structures}, Annals of Mathematics Studies, Vol. 169, Princeton University Press, 2009.
\bibitem[Sab21]{Sab21} Claude Sabbah, \textit{Harmonic Metrics and Wild Ramification}, Annales scientifiques de l'ENS, 54(5), 1235–1294, 2021.
\bibitem[Moc16]{Moc16} Takurō Mochizuki, \textit{Kobayashi-Hitchin Correspondence for Twisted Harmonic Bundles}, Journal of Differential Geometry, 102(2), 273–397, 2016.
\bibitem[Bar16]{Bar16} Clark Barwick, On the algebraic K-theory of higher categories, \textit{Journal of Topology}, 9(1), 245–347, 2016.
\bibitem[BGT13]{BGT13} Anthony Blanc, David Gepner and Gonçalo Tabuada, A Universal Characterization of Higher Algebraic K-Theory, \textit{Annals of Mathematics}, 178(3), 1147–1200., 2013.
\bibitem[HLL18]{HLL18} Tobias Harrer, Daniel Lucas and Simon Pepin Lehalleur, On the Nori Motivic Groupoid, \textit{Compositio Mathematica}, 154(9), 1827–1877, 2018.
\bibitem[Lur17]{Lur17} Jacob Lurie, \textit{Higher Algebra}, 2017. Available at https://www.math.ias.edu/~lurie/papers/HA.pdf.
\bibitem[Lur18]{Lur18} \underline{\hspace{1cm}}, \textit{Spectral Algebraic Geometry}, 2018, Institute for Advanced Study. Available at https://www.math.ias.edu/~lurie/papers/SAG-rootfile.pdf.
\bibitem[Dre15]{Dre15} Lewis Drew, \textit{The Hodge $\infty$-category and motivic homotopy theory}, Preprint, 2015. Available at https://drew-lewis.com/publications/.
\bibitem[IKMZ16]{IKMZ16} Ilia Itenberg, Ludmil Katzarkov, Grigory Mikhalkin and Ilia Zharkov, Tropical homology. arXiv preprint arXiv:1607.06396, 2016.
\bibitem[Voi02]{Voi02} Claire Voisin, \textit{Hodge Theory and Complex Algebraic Geometry I}, Cambridge Studies in Advanced Mathematics 76, Cambridge University Press, 2002. 
\bibitem[AMRT75]{AMRT75} Avner Ash, David Mumford, Michael Rapoport and Yung-Sheng Tai, \textit{Smooth compactifications of locally symmetric varieties}, Math. Sci. Press, Brookline, MA, 1975.
\bibitem[And96]{And96} Yves Andr\'e, \textit{Pour une théorie inconditionnelle des motifs}, Publications Mathématiques de l'Institut des Hautes Études Scientifiques, 83 (1996), 5–49.
\end{thebibliography}
\end{document}